\documentclass[12pt]{amsart}

\usepackage{amsmath}
\usepackage{amssymb}
\usepackage{graphicx}
\usepackage{pifont}
\usepackage{epsfig}
\usepackage{verbatim}
\usepackage[english]{babel}
\usepackage{fancyhdr, blindtext}
\usepackage{float}
\usepackage{subfig}

\usepackage{graphics}
\usepackage{graphicx,epsfig}
\usepackage{epsfig}
\usepackage{graphicx}

\usepackage{a4wide}
\usepackage{amsmath}
\usepackage{amsfonts}

\usepackage{enumerate}
\usepackage{authblk}
\usepackage{float}
\usepackage{subfig}
\usepackage{amsaddr}

\graphicspath{{./figures_mathbio/}}

\begin{document}

\newcommand\bes{\begin{eqnarray}}
\newcommand\ees{\end{eqnarray}}
\newcommand\bess{\begin{eqnarray*}}
\newcommand\eess{\end{eqnarray*}}
\newcommand{\ve}{\varepsilon}
\newtheorem{definition}{Definition}
\newtheorem{theorem}{Theorem}[section]
\newtheorem{lemma}{Lemma}[section]
\newtheorem{proposition}{Proposition}[section]
\newtheorem{remark}{Remark}[section]
\newtheorem{corollary}{Corollary}[section]
\newtheorem{example}{Example}[section]
\title[Complex oscillatory patterns]
{\bf Complex oscillatory patterns near singular Hopf bifurcation in a two time-scale ecosystem}

\author{Susmita Sadhu}
\address{Department of Mathematics,
Georgia College \& State University, Milledgeville, GA 31061, USA}
\email{susmita.sadhu@gcsu.edu}
\urladdr{http://faculty.gcsu.edu/custom-website/susmita-sadhu}


\thispagestyle{empty}


\begin{abstract}
\noindent   
We consider an ecological model consisting of two species of predators competing for their common prey with explicit interference competition.   With a proper rescaling, the model is portrayed as a singularly perturbed system with one-fast  (prey dynamics) and two-slow variables (dynamics of the predators).   The model exhibits variety  of rich and interesting dynamics, including, but not limited to mixed mode oscillations (MMOs), featuring  concatenation of small and large amplitude oscillations, relaxation oscillations and bistability  between a semi-trivial equilibrium state and a coexistence oscillatory state. Existence of co-dimenison two bifurcations such as fold-Hopf and generalized Hopf bifurcations make the system further intriguing. More interestingly, in a neighborhood of {\emph{singular Hopf}} bifurcation, long lasting transient dynamics in  form of chaotic MMOs or relaxation oscillations  are observed as the system approaches the periodic attractor born out of supercritical Hopf bifurcation or a semi-trivial equilibrium state respectively. The transient dynamics could  persist for hundreds or thousands of generations before the ecosystem experiences a regime shift. The time series of population cycles with different types of irregular oscillations arising in this model stem from a biological realistic feature, namely, by the variation in the intraspecific competition amongst the predators. To explain these oscillations, we use bifurcation analysis and methods from {\emph{geometric singular perturbation theory}}.

\end{abstract}

\maketitle

Key Words. Mixed-mode oscillations, singular Hopf bifurcation,  long transients, basins of attraction, fold-Hopf bifurcation, chaotic dynamics.

\vspace{0.15in}



\section{Introduction}

Variability in population abundance is one of the most ubiquitous phenomena in natural populations. These dynamics typically involve rapid transitions from one oscillatory state to another, featuring temporal patterns of small fluctuations interspersed with large amplitude oscillations \cite{KC, Peet}, as observed in population cycles of Pacific sardine and northern anchovy \cite{Baumgartner}, 
 forest insects such as larch bud moth, gypsy moth and cankerworms \cite{AC,KBJP}, agricultural pests such as desert locusts \cite{C, W1} and so forth. Using singular perturbation theory,  regular cycles of population outbreaks and collapses are modeled by relaxation oscillation cycles, commonly known as boom and bust cycles \cite{LXY, MR, RM}.  However, since  events of population outbreaks are not regular, vary in intensities and generally not periodic in time, they are perhaps not best represented by relaxation oscillations.  A more realistic representation of these cycles can be mixed-mode oscillations (MMOs), which are complex oscillatory patterns consisting of one or more small amplitude oscillations (SAOs) followed by large amplitude oscillations (LAOs) which are of relaxation type \cite{DGKKOW}.  MMOs have been widely studied in various physical systems, including, but not limited to computational neuroscience models, autocatalytic and  electrochemical reactions (see \cite{BKR, CR, K, K11, LRV, MS, RWK} and the references therein for  examples of mathematical models exhibiting MMOs). In ecology, there have been relatively few studies on such complex  oscillatory patterns \cite{BK, KC, Peet, S1, SCT, SK}.
 
Traditionally, theoretical ecological studies have focussed on asymptotic analysis to study population dynamics. However, there has been an increasing recognition of  importance of understanding the properties of transient dynamics in complex ecosystems  \cite{Hastings, morozovetal}. One of the major issues in ecology is predicting  sudden changes in population abundance, commonly referred to as regime shifts   \cite{SC,Scheffer},  as  they often result in population collapse and extinctions.    Long term transient dynamics on ecological timescales can offer novel perspectives in underlying mechanisms behind regime shifts in an apparently constant environment. The present work unifies the two phenomena that appear in natural populations, namely MMOs featuring oscillations on different timescales, and their occurrence as long transient dynamics over ecological timescales.

In this paper, we study the dynamics between three interacting species, namely two predators competing for their common prey with explicit Lotka-Volterra type interference competition between the predators. 
 Taking the ratios between  birth rates of the predators to the prey extremely small, separation of time scales is introduced into the system as a singular parameter. This then brings the  model in the framework of singularly perturbed system of equations with one-fast and two-slow variables, whereby the prey exhibits fast dynamics and the predators have slow dynamics. Similar assumptions have been made in many ecological models, including but not limited to forest pest models 
\cite{BK}, tri-trophic food chain models \cite{BD1, DH}, age-structured predator-prey models with dormancy of predators \cite{KC} and the well-known Rosenzweig-MacArthur predator-prey model \cite{PAGK}.  Though the model studied in this paper is a generic ecological model, it can give rise to a host of new and interesting patterns of behavior. 
These dynamics seem to represent natural population cycles. 
A primary goal of this paper is to study various kinds of complex oscillatory dynamics that emerges in this model and use tools from geometric singular perturbation theory and bifurcation analysis to explain them.


Treating the intraspecific competition amongst one of the species of predators as  the main control parameter and the predation efficiency of the other species of  predators as the secondary parameter  ($\beta_1$), the local and global bifurcation structures in the system  are explored. We find different dynamical regimes on which interesting ecological phenomena occur, and compare the bifurcation structures of the model by considering different values of the predation efficiency of the other predator (denoted by $\beta_2$). Bifurcations such as saddle-node, Hopf, period-doubling, torus, fold-Hopf and generalized Hopf occur in the system and contribute to unfolding different interesting phenomena exhibited by the model.
 
 One of the interesting dynamics observed in this model is the existence of long term transients in form of chaotic MMOs consisting of very long epochs of SAOs near the {\emph{singular-Hopf}} point, also referred to as {\emph{folded saddle-node of type II}} (FSN II) bifurcation point \cite{KPK, KW}. The  transient MMO dynamics persists for thousands of generations and eventually asymptotes to the small amplitude periodic attractor born out of a supercritical Hopf bifurcation. The stable periodic attractor possesses a large basin of attraction and interferes with the global return mechanism that is responsible for the relaxation dynamics of the MMOs. The SAOs associated with the chaotic transients are organized by a slow passage near special points in phase space  such as {\emph{folded node singularities}} \cite{KW, SW} and are further influenced by the local vector field around a saddle-focus equilibrium which lies in a close vicinity to the folded node singularity. The existence of long lasting transient dynamics reflects that the dynamics on an ecological timescale can be totally different than asymptotic dynamics. 
 The occurrence of  a supercritical Hopf bifurcation, which is regarded as a {\emph{safe bifurcation}} in context of ecological regime shifts \cite{SC,Scheffer}  can be misleading on an ecological timescale.

In another parameter regime, where the coexistence equilibrium undergoes a subcritical Hopf bifurcation, transient dynamics in form of relaxation oscillations are observed before the system approaches a boundary equilibrium state, where one of the predators becomes extinct. In a nearby neighborhood of the subcritical Hopf bifurcation, a torus bifurcation occurs, giving birth to a  family of stable periodic orbits. In this regime, the system exhibits a bistable behavior, namely depending on the initial densities, the model  could either result into coexistence of species in an oscillatory state or lead to a competitive exclusion of one of the species of  predators.

Long term transients  in form of  MMOs or relaxation oscillations shed light on the underlying inherent uncertainties associated with population dynamics. 
These can be ecologically alarming, as abrupt large population fluctuations may lead to outbreaks of harmful pests such as desert locusts or forest insects or bring a species perilously close to extinction under a seemingly constant environment.  In fact, such long lasting transient dynamics can provide an alternate explanation of a sudden transition from one state to another in an ecosystem which may show an apparently stable dynamics or long statis \cite{Hastings}.   
 In this model, the local dynamics near the equilibrium can be used to identify early warning signs of a  large fluctuation which has been explored in a companion paper \cite{Sadhunew}.   

In multiple timescale systems with one fast and two or more slow variables, the two main mechanisms attributed to the generation of MMOs are the {\emph{generalized canard phenomenon}} \cite{BKW, SW}  and a {\emph{singular Hopf bifurcation}} \cite{DGKKOW, G}. In \cite{G}, the differences in the characteristics of SAOs in an MMO orbit  due to singular Hopf bifurcation and those created due to folded node singularity is discussed: namely, the small oscillations associated with singular Hopf start with low amplitude and then grow before entering into relaxation regime, whereas the SAOs associated with folded nodes decrease and then increase in amplitudes. However in this model, MMOs with both characteristics are seen near FSN II bifurcation, where a folded node lies in the vicinity of the equilibrium. Similar dynamics have been previously observed in neuroscience models \cite{CR, RWK}.  More interestingly in an another parameter regime, MMOs with another unique characteristic is observed:  namely, the durations of the episodes of small oscillations between two large oscillations vary by orders of magnitude, coupled with large variations in the amplitudes of the small oscillations. A major focus on this paper is to study different MMO patterns that arise in this model near FSN II bifurcation.  
We use the presence of different timescales in the system to analyze the mechanism for the complex oscillations. The analyses performed at the singular limit are  used to explain the behaviors seen away from the limit.

The paper is organized as follows. In Section 2, we introduce and properly scale the model. The assumptions and physical significance of each parameter are also discussed.   Section 3 focusses on examples of different MMO patterns exhibited by the model either as asymptotic dynamics or as long term transients. In parameter regimes, where MMOs occur as prolonged transients,  ``transient basins of attraction" are computed and compared.  Relaxation oscillations as long term transients in a different parameter regime are also studied in Section 3. Section 4 focusses on bifurcation analysis of the model, where interesting dynamics near the folded saddle-node singularity of type II are explored. A geometric singular perturbation approach and relevant background review is given in Section 5, which provides us a framework of local mechanisms that are responsible for MMOs.  The numerical results are explained using this framework. 
Finally, we discuss our results and summarize our conclusions in the last section of the paper.


\section{The Model}

The model studied in this paper is  an extension of the model considered in \cite{RARD} with additional direct interactions between the predators. The model reads as follows:
\begin{eqnarray}\label{maineq} \left\{ 
\begin{array}{ll}\label{1}
      \frac{dX}{dT} &= rX\left(1-\frac{X}{K}\right)-\frac{p_1XY}{H_1+X}-\frac{p_2XZ}{H_2+X}\\
        \frac{dY}{dT} &= \frac{b_1p_1XY}{H_1+X}-d_1Y -a_{12}YZ\\
         \frac{dZ}{dT} &= \frac{b_2p_2XZ}{H_2+X}-d_2Z -a_{21}YZ -mZ^2
       \end{array} 
\right. \end{eqnarray}
under the initial conditions
\bes \label{maineqic}  X(0)=\tilde{X}\geq 0,\ Y(0)=\tilde{Y}\geq 0, \ Z(0)=\tilde{Z}\geq 0,
\ees
where $X$ represents the population density of the prey and $Y$, $Z$ represent the densities of the two species of  predators. The parameters $r$ and $K$ represent  the intrinsic growth rate and the carrying capacity of the prey, $p_1$ is the maximum per-capita predation rate of $Y$, $H_1$ is the semi-saturation constant which represents the prey density at which $Y$ reaches half of its maximum predation rate ($p_1/2$),  $b_1$ and $d_1$ are the birth-to-consumption ratio and  per-capita natural death rate of $Y$ respectively. The parameter $a_{12}$ is the rate of adverse effect of $Z$ on $Y$, which we will refer to as the interspecific competition rate. The other parameters $p_2, b_2, d_2, H_2, a_{21}$ are defined analogously for $Z$. We assume that $Z$ is territorial and/or is a  social predator, and therefore experiences more competition for space and resources in comparison to $Y$ (which is assumed to be  non-territorial or is a solitary predator). The term $mZ^2$ in the $Z$ equation accounts for this intraspecific competition, which may include lethal fighting and cannibalism, and measures the density dependent mortality rate in the class of $Z$. Examples of species governed by (\ref{maineq}) may include small mammals such as  rabbits/rodents preyed upon by raptors such as eagles ($Y$) and larger mammals such as wolves ($Z$) or planktonic and small nektonic species such as fish and squids ($X$) preyed upon by sharks ($Y$) and dolphins ($Z$). 
      
The density dependent mortality term in system (\ref{maineq}) is the defining term  that distinguishes the possibility of existence of a positive equilibrium state from competitive exclusion of species.  In case of two competitors competing for a common species,  the stronger competitor can outcompete the other; a phenomenon commonly referred to as the principle of competitive exclusion \cite{H}. Such dynamics are  known to occur in Rosenzweig-MacArthur and Armstrong and McGhee predator-prey models \cite{HHW, MA, MR}. In these models, the only coexistence state that the systems admit are in the form of large oscillation cycles \cite{LXY}.  However, a density dependent mortality term can prevent competitive exclusion and ensure  long-term survival of the species besides guaranteeing  the existence of its competitor \cite{RARD}. Thus, the term $mZ^2$ is ecologically significant and we will treat the coefficient of intraspecific competition as the key parameter in this work.

With the following change of variables and parameters:
\begin{eqnarray} {\nonumber} t&=&rT, \  x= \frac{X}{K}, \ y=\frac{p_1Y}{rK},\ z=\frac{p_2Z}{rK}, \ \zeta_1 = \frac{b_1p_1}{r},  \ \zeta_2 = \frac{b_2p_2}{r}, \ \beta_1 =  \frac{H_1}{K}, \\
{\nonumber}   \beta_2 &=&\frac{H_2}{K}, \ c = \frac{d_1}{b_1p_1},\  d= \frac{d_2}{b_2p_2}, \   h= \frac{mZ_0}{b_2p_2}, \ \alpha_{12} = \frac{a_{12}Z_0}{b_1p_1},\  \alpha_{21} = \frac{a_{21}Y_0}{b_2 p_2},
\end{eqnarray}
where
\begin{eqnarray} {\nonumber} Y_0=\frac{rK}{p_1},\quad  Z_0=\frac{rK}{p_2},
\end{eqnarray}
system $(\ref{maineq})$ takes the following  dimensionless form:
\begin{eqnarray}\label{nondim1}    \left\{
\begin{array}{ll} {x'}&= x\left(1-x-\frac{y}{\beta_1+x}-\frac{z}{\beta_2+x}\right)\\
 {y'}&=\zeta_1 y\left(\frac{x}{\beta_1+x}-c-  \alpha_{12} z \right)\\
   {z'} &= \zeta_2  z\left(\frac{x}{\beta_2+x}-d -\alpha_{21} y-hz \right),
       \end{array} 
\right. \end{eqnarray}
where the primes denote differentiation with respect to the time variable $t$. Similar scaling variables were considered in a three-trophic food chain model by Deng \cite{BD1}.  
The quantity $Y_0$  will be interpreted as the maximum predation capacity of  $Y$ \cite{BD1}.  More precisely, at this capacity,  $Y$ will consume $rK$ number of prey (which measures the  reproduction capacity of the prey population)  at the maximum predation rate $p_1$.  Similarly, $Z_0$ will be considered as the maximum predation capacity of $Z$. We will assume the following conditions on the parameters: 

(A) The maximum per capita birth rate of the prey is much higher than the per capita birth rates of the predators, i.e. $b_1p_1 << r$ and $b_2p_2<<r$, thus yielding  $0<\zeta_1, \zeta_2<<1$. This is usually observed in many ecosystems \cite{MR1}. 
Examples include insects and their avian predators, rodents and their aerial or ground-based predators and so forth. For simplicity,  we will assume in our model that $\zeta_1= \zeta_2=\zeta$ (say). Similar dynamics are obtained when $\zeta_1\neq \zeta_2$, but have the same order.

(B) The parameters $c$ and $d$ satisfy the inequality $0<c,\ d<1$, which implies that the growth rates of the predators are greater than their death rates. This is a default assumption otherwise the predators would die out faster than they could reproduce even at their maximum reproduction rate.

(C) The parameters $\beta_1$ and $\beta_2$ are dimensionless semi-saturation constants measured against the prey's carrying capacity. We will assume that both predators are efficient, and hence they will reach the half of their maximum predation rates before the prey population reaches its carrying capacity, thus yielding $0< \beta_1,\beta_2<1$.

(D) The parameters $\alpha_{12}$ and $\alpha_{21}$ are dimensionless interspecific competition coefficients measuring the interference effect of $Z$ on $Y$ and of $Y$ on $Z$ respectively.  We will assume that the maximum interference effect of one species on the other does not exceed the intrinsic growth rate of the other species, yielding $0<\alpha_{12},\alpha_{21}<1$.

 
Under the assumptions (A)-(D), system $(\ref{nondim1})$  transforms to a singular perturbed system of equations with two time scales, where the prey exhibits fast dynamics and the predators exhibit slow dynamics. 

We rewrite system $(\ref{nondim1})$  as
\begin{eqnarray}\label{nondim2}    \left\{
\begin{array}{ll}  {x'} &= x\left(1-x-\frac{y}{\beta_1+x}-\frac{z}{\beta_2+x}\right):=x\phi(x,y,z)\\
    {y'}&= \zeta y\left(\frac{x}{\beta_1+x}-c-\alpha_{12} z\right):=\zeta y\chi(x,z)\\
    {z'}&=  \zeta z\left(\frac{x}{\beta_2+x}-d - \alpha_{21} y -hz  \right):= \zeta z\psi(x,y,z),
      \end{array} 
      \right. 
\end{eqnarray}
where $\phi=0$, $\chi=0$, and $\psi=0$ are the nontrivial $x$, $y$, and $z$-nullclines respectively. On rescaling $t$ by $\zeta$  and letting $s=\zeta t$, system $(\ref{nondim2})$ can be reformulated as

\begin{eqnarray}\label{nondim3}    \left\{
\begin{array}{ll} \zeta \dot{x} &=x\phi(x,y,z)\\
    \dot{y}&=y\chi(x,z)\\
    \dot{z}&=  z\psi(x,y,z),
       \end{array} 
\right. 
\end{eqnarray}
where the overdot denotes differentiation with respect to the variable $s$. The variables $t$ and $s$ are referred to as the fast and slow time variables respectively. The parameter $\zeta$ can be regarded as the separation of time scales.

\begin{figure}[h!]     
  \centering   
      {\includegraphics[width=7.50cm]{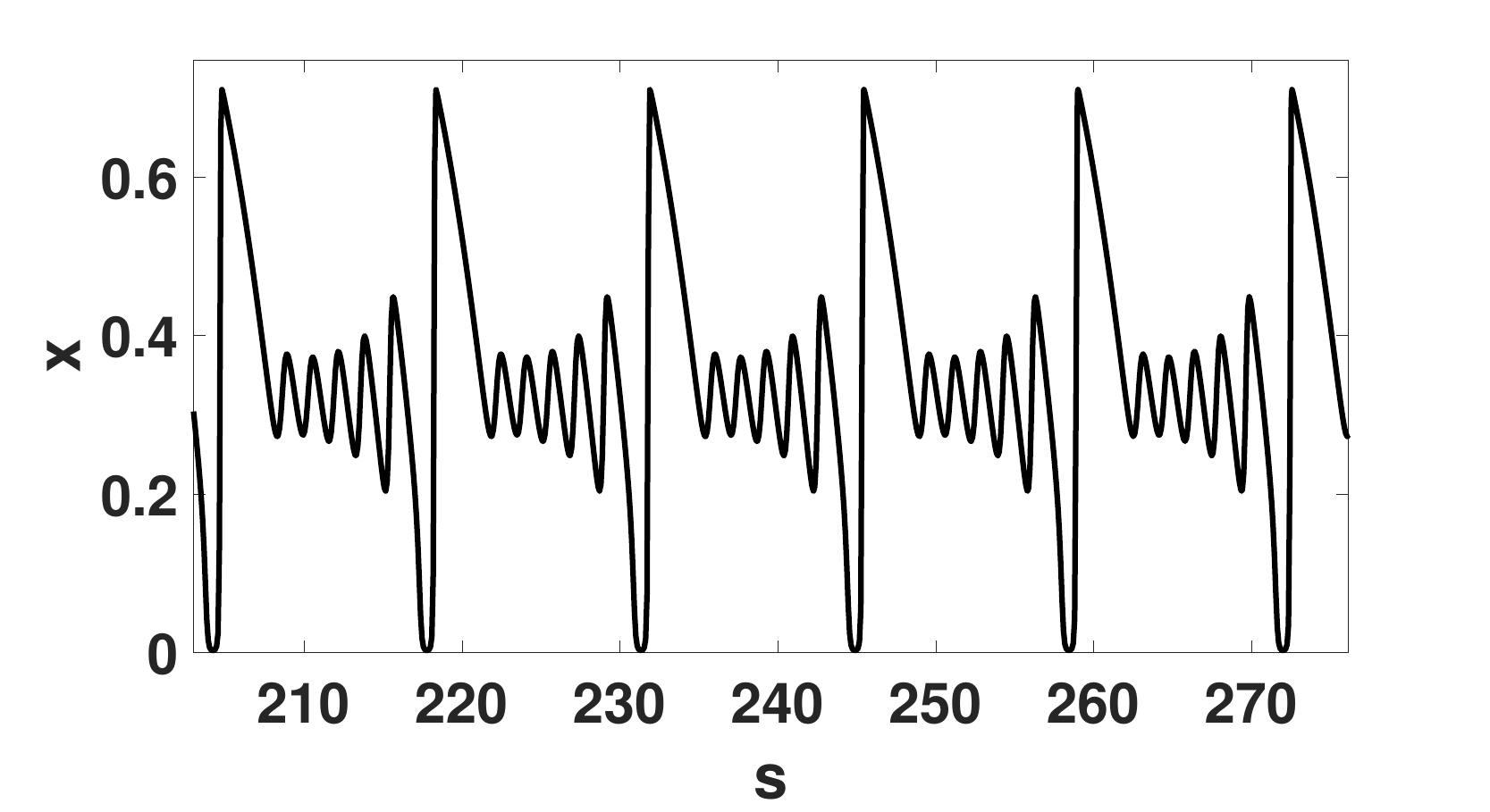}}\qquad
      {\includegraphics[width=7.50cm]{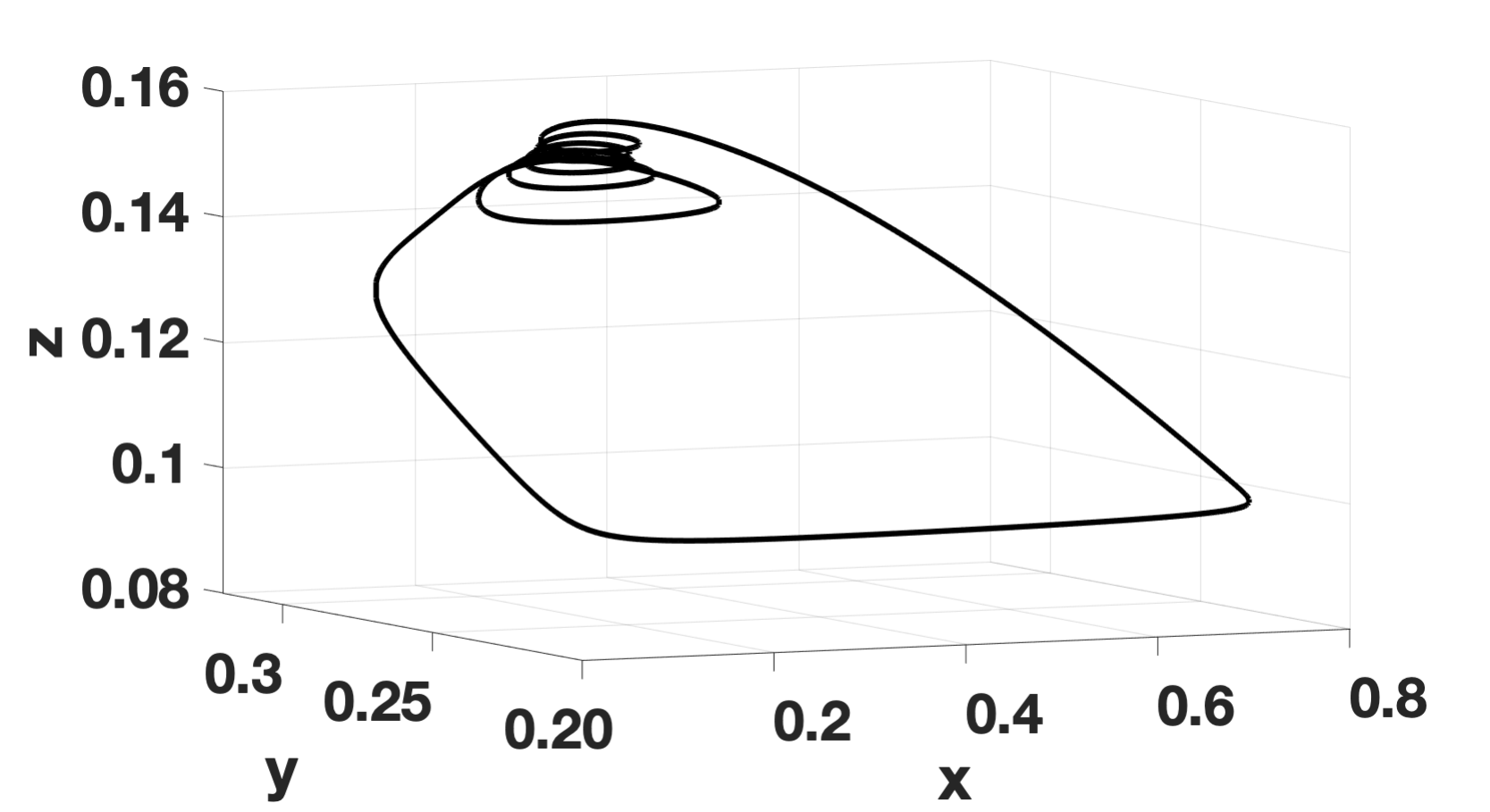}}
  \caption{Time series and phase portrait of a $1^5$  MMO pattern observed in system (\ref{nondim3}) for  
  $\beta_1=0.35$,  $h=1.65$ and the other parameters are as in (\ref{parvalues}). 
  }
  \label{timeseries_example}
\end{figure}

Throughout the paper, we will fix the parameter values to
\bes \label{parvalues}
\zeta =0.01, \ \beta_2 =0.35,\ c=0.4,\ d=0.21,\ \alpha_{12} =0.5, \ \alpha_{21}=0.1,
\ees
or to
\bes \label{parvalues1}
\zeta =0.01, \ \beta_2 =0.6,\ c=0.4,\ d=0.21,\ \alpha_{12} =0.5, \ \alpha_{21}=0.1,
\ees
and treat the intraspecific competition coefficient $h$ as the primary control parameter, and the predation efficiency $\beta_1$ as the secondary parameter.
The parameter values chosen here are for illustrative purposes to demonstrate the complex oscillatory patterns that could occur in a complex ecosystem as shown in figure \ref{timeseries_example}. 

System (\ref{nondim3}) has $8$ free parameters which make the qualitative analysis very challenging. In the absence of interference competition, i.e. $\alpha_{12}=\alpha_{21}=0$, there exists a unique coexistence equilibrium in a suitable parameter regime. MMO orbits and chaotic dynamics exist as a control parameter is varied. This case has been studied in details in \cite{SCT}. 
 In this paper, we investigate the dynamics  in a neighborhood of the {\emph{singular Hopf}} point for $\alpha_{12}, \alpha_{21} \neq 0$. To this end, we take all the parameter values to be nonzero such that  the intersection of the non-trivial nullclines $\phi=0, \chi=0$ and $\psi=0$ produces 
 equilibria that lie in the positive octant. These equilibria will be referred to as coexistence or non-trivial equilibria and will be denoted by $E_i^*$.
\section{Complex oscillatory dynamics in system (\ref{nondim3})}
\subsection{MMOs as asymptotic dynamics}

\begin{figure}[h!]     
  \centering 
  {\includegraphics[width=9.5cm]{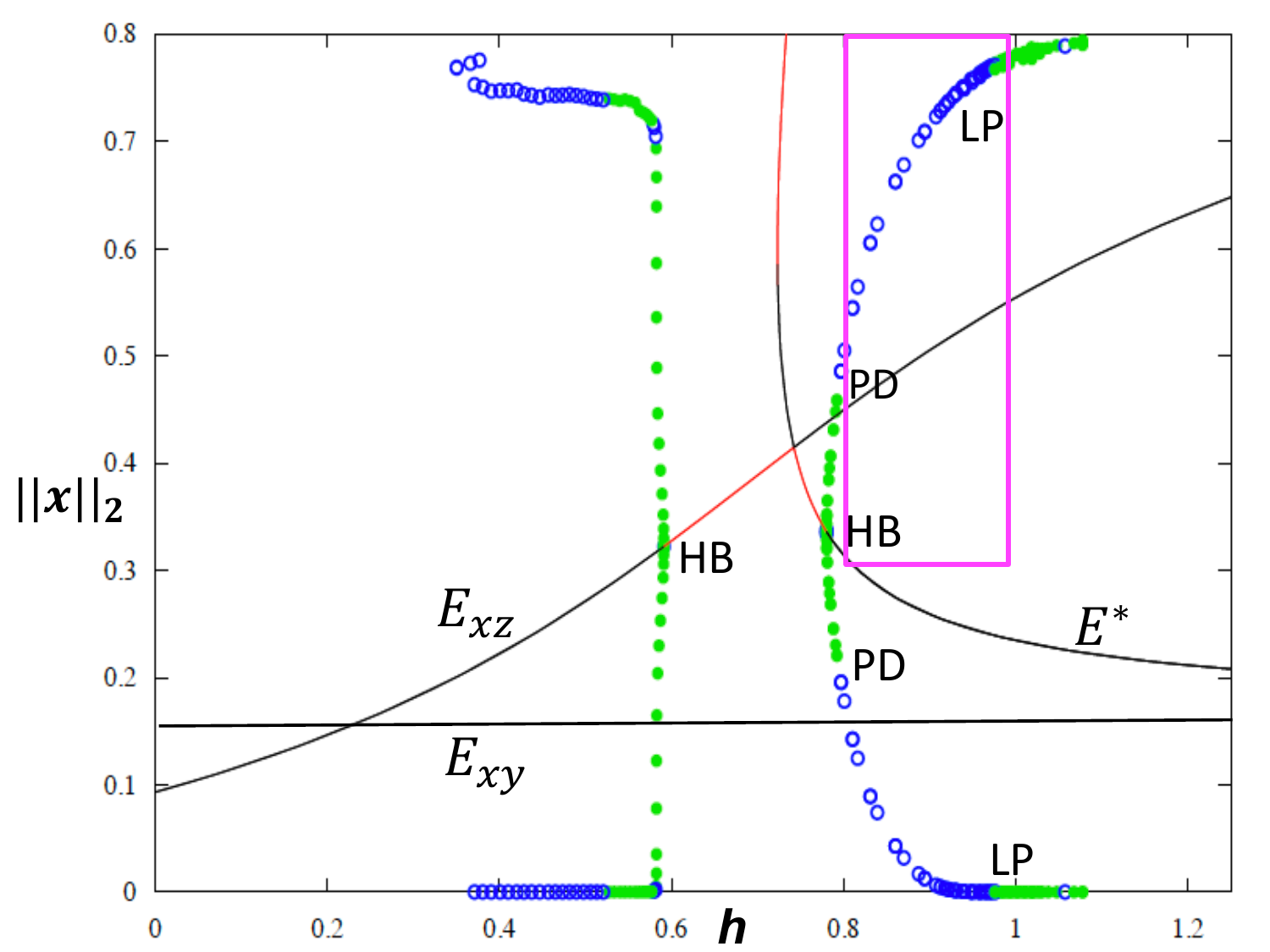}}
  \caption{A one-parameter bifurcation diagram in $h$ for $\beta_1=0.25$ and other parameters as in (\ref{parvalues}). The red (black) represent stable (unstable) branches of equilibria, where $E^*$ denotes the coexistence equilibrium, $E_{xy}$ and $E_{xz}$ are the boundary equilibria on the $xy$ and $xz$ planes respectively. The green (blue) circles represent stable (unstable) limit cycles. HB: Hopf bifurcation, PD: period doubling, LP: saddle node bifurcation of limit cycles. See the next figure for a detailed view of  the inset. The plot was generated in XPPAUT \cite{E}.}
  \label{one_par_xpp}
\end{figure}

Fixing $\beta_1=0.25$  and the other parameter values as in  (\ref{parvalues}), we consider the bifurcations that arise by varying the intraspecific competition $h$. Using AUTO, a one-parameter bifurcation diagram was calculated as shown in figure (\ref{one_par_xpp}), where $h$ is the continuation parameter and $L_2$ norm refers to the standard Euclidean norm.  We note from the bifurcation diagram in figure \ref{one_par_xpp} that the coexistence equilibrium point $E_1^*$ exists for $h>0.7425$. In the singular limit of system (\ref{nondim3}),  FSN II bifurcation occurs at $h\approx 0.7785$ (also see figure \ref{one_par_desing}),  
where $E_1^*$ is still an attractor of the full system.  A supercritical Hopf bifurcation occurs at $h \approx 0.7803$, giving birth to a family of stable periodic orbits $\Gamma_h$.  The equilibrium $E_1^*$ is now a saddle-focus with one negative and two complex (with positive real parts) eigenvalues. The periodic orbits $\Gamma_{h}$ grow in size and undergo a cascade of period-doubling bifurcations for $h \in (0.797, 0.798)$, resulting into small amplitude chaotic invariant sets.   

  \begin{figure}[h!]     
  \centering 
{\includegraphics[width=8.5cm]{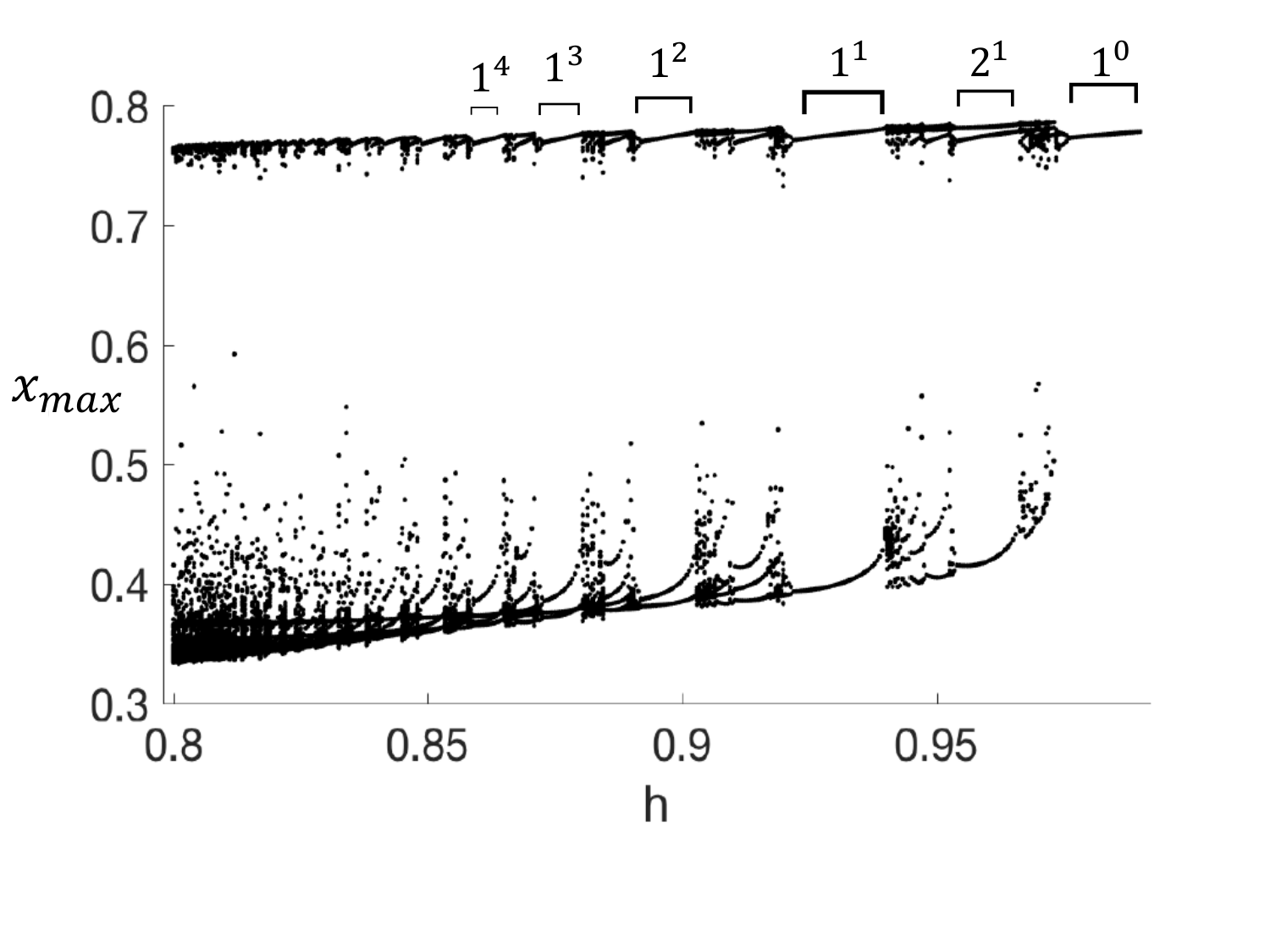}}\quad
 \caption{A detailed one-parameter bifurcation structure of MMOs corresponding to the inset in figure \ref{one_par_xpp}.  The maximum value of $x$ is shown in the $y$-axis.}
 \label{full_bif_025}
\end{figure}

\begin{figure}[h!]     
  \centering   
  \subfloat[MMO orbit]{\includegraphics[width=7.05cm]{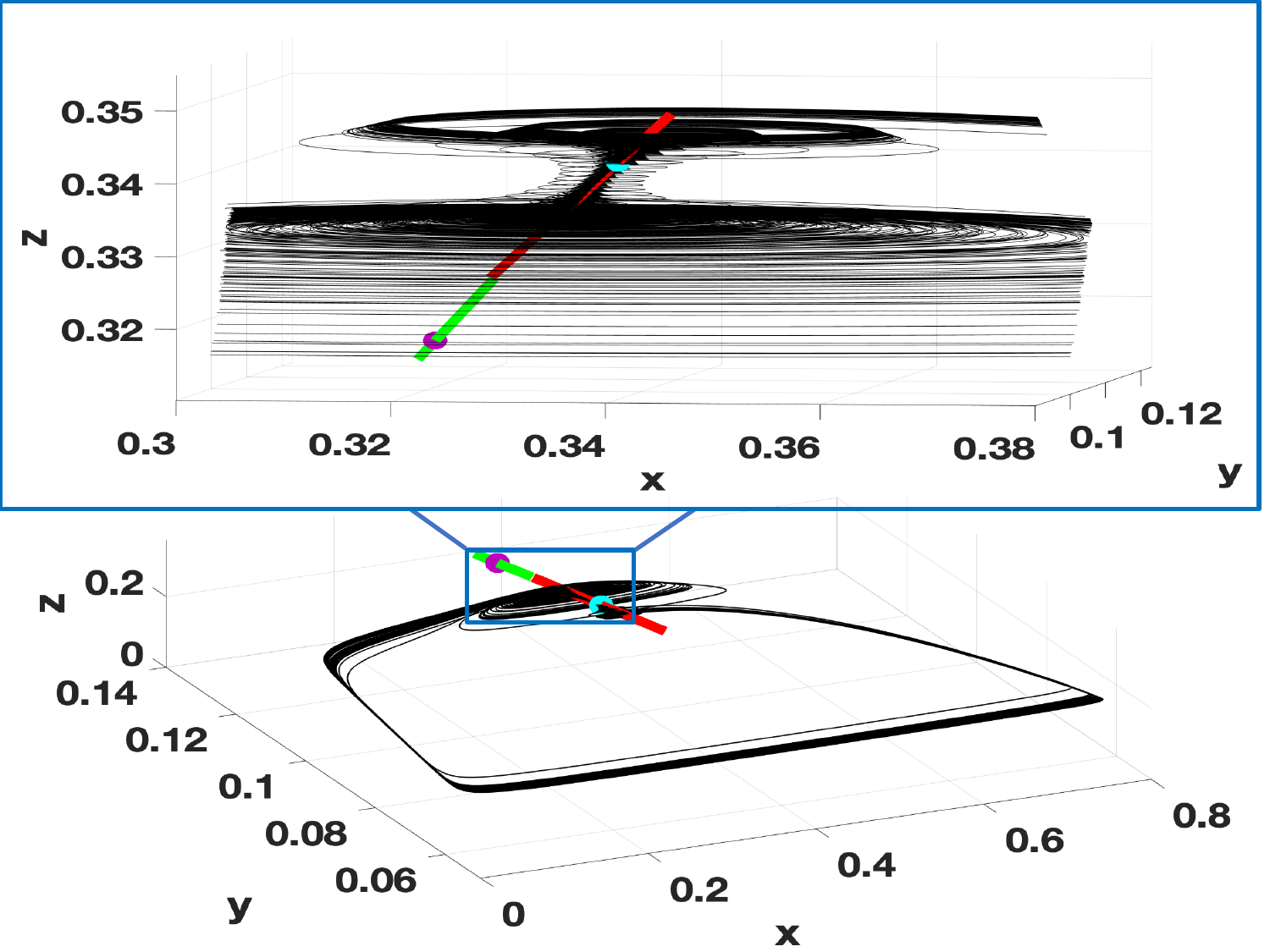}}\qquad
      \subfloat[Time series in $x$]{\includegraphics[width=7.0cm]{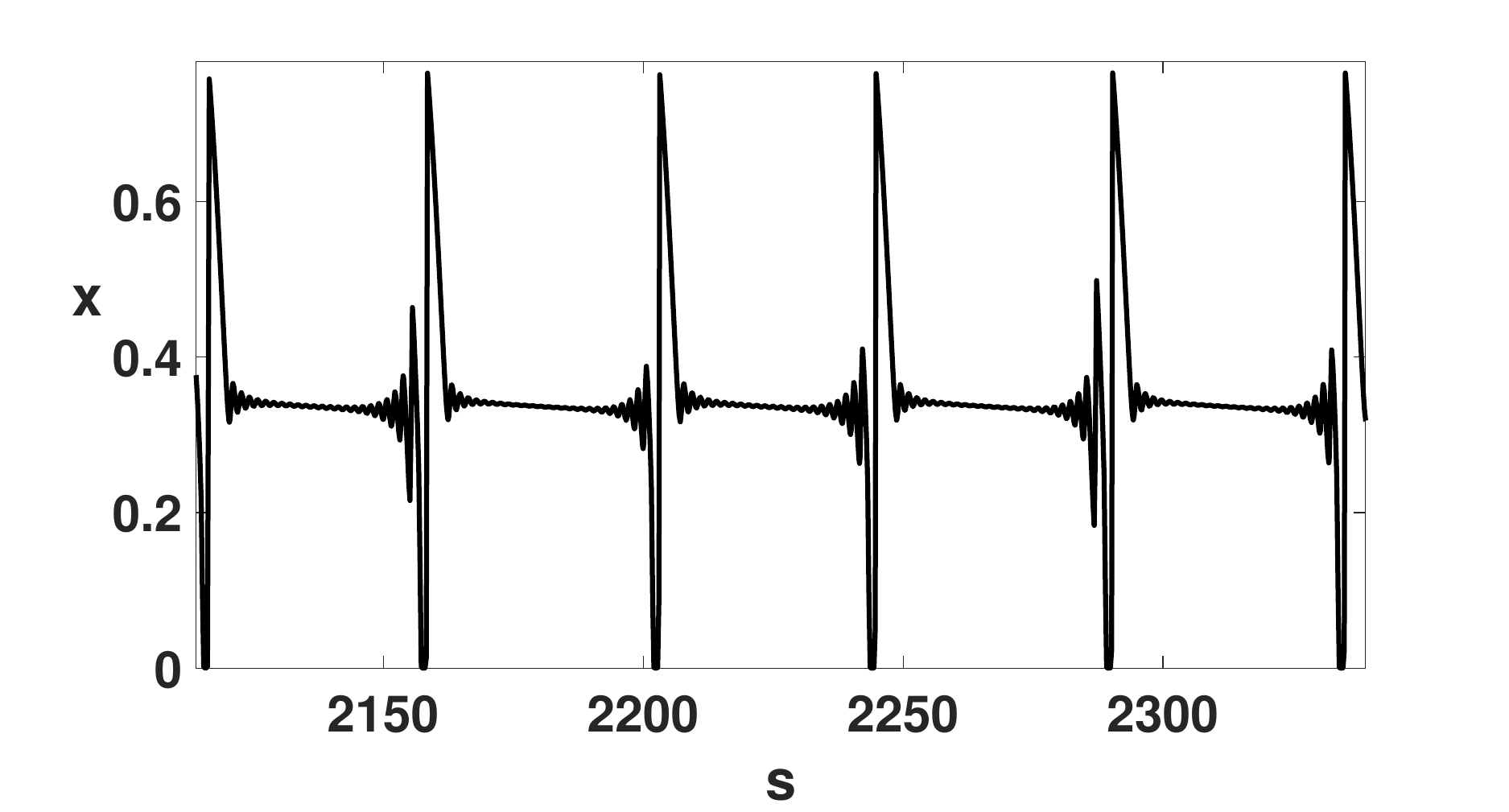}}
  \caption{ (A): Phase portrait of a chaotic MMO orbit. The inset contains a zoomed view of the orbit near the folded node (cyan dot) as it rotates around the weak eigendirection (red) while getting pushed towards the equilibrium (magenta) along its stable eigendirection (green). (B): Corresponding time series in $x$. 
  Here $\beta_1=0.25$,  $h=0.8$ and the other parameters are as in (\ref{parvalues}). 
  }
  \label{timeseries_08_819}
\end{figure}

 Past the regime of period doubling bifurcations of $\Gamma_h$, chaotic dynamics featuring small amplitude and large amplitude oscillations are observed.  On further increasing $h$, a sequence of period doubling and saddle-node bifurcations of limit cycles produce periodic MMO orbits predominantly with signatures $1^s$ where $s \in \mathbb{N}$ can be very large, and other complicated signatures, which are usually a mix of $1^i$ and $1^{i+1}$ as shown in figure \ref{timeseries_08_819}.  The SAOs in such an MMO orbit are organized by the rotational properties of the {\emph{weak canard}} as well as by the unstable manifold of $E_1^*$.  We will  discuss the local mechanisms generating these SAOs in Section 5. Relaxation oscillations are observed for $h>0.98$.  A detailed view of the bifurcation structure of the MMOs is given in  figure \ref{full_bif_025}.  
 \begin{figure}[h!]     
  \centering   
              \subfloat[]{\includegraphics[width=7.5cm]{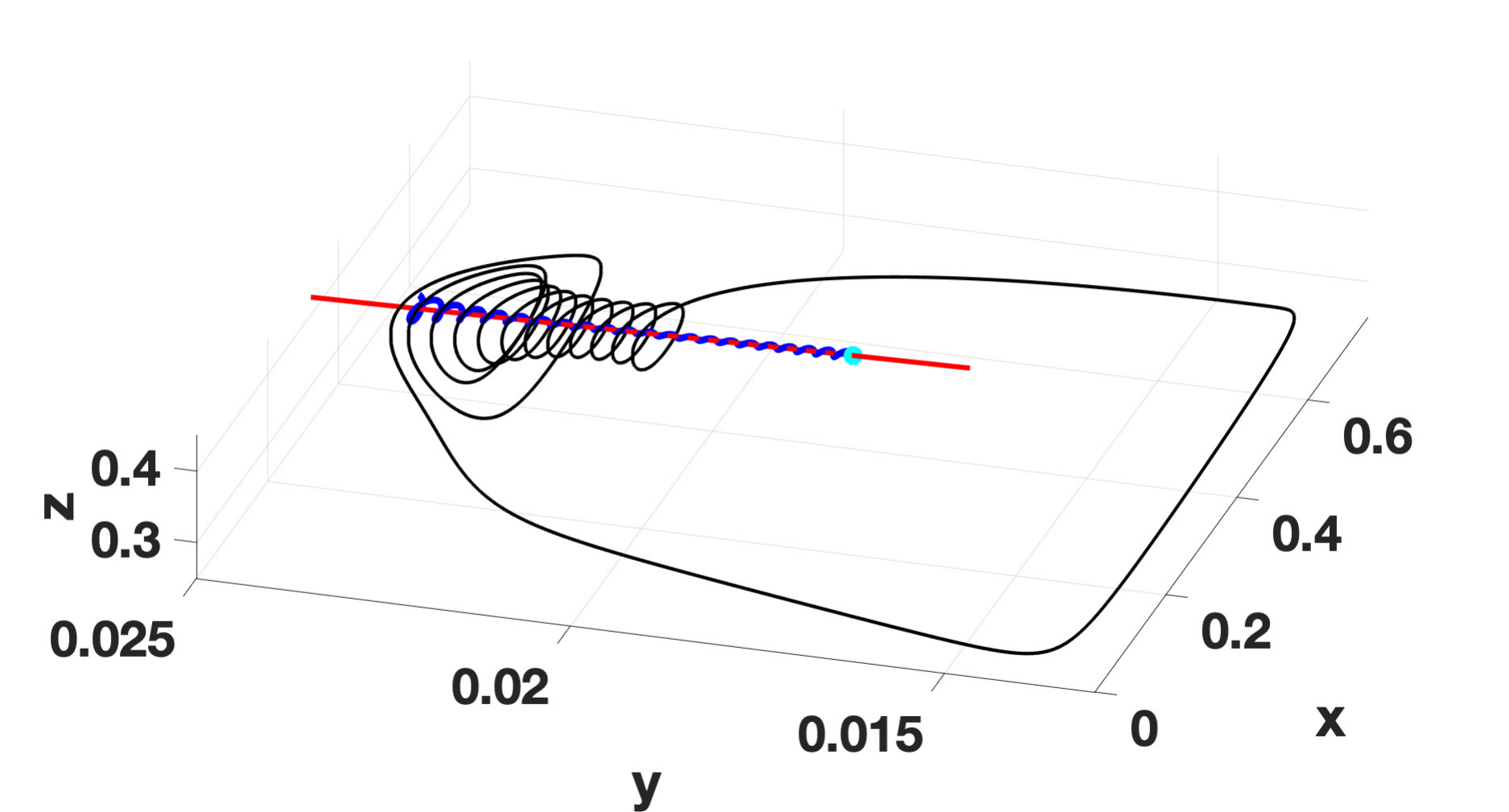}}\qquad
          \subfloat[]{\includegraphics[width=5.7cm]{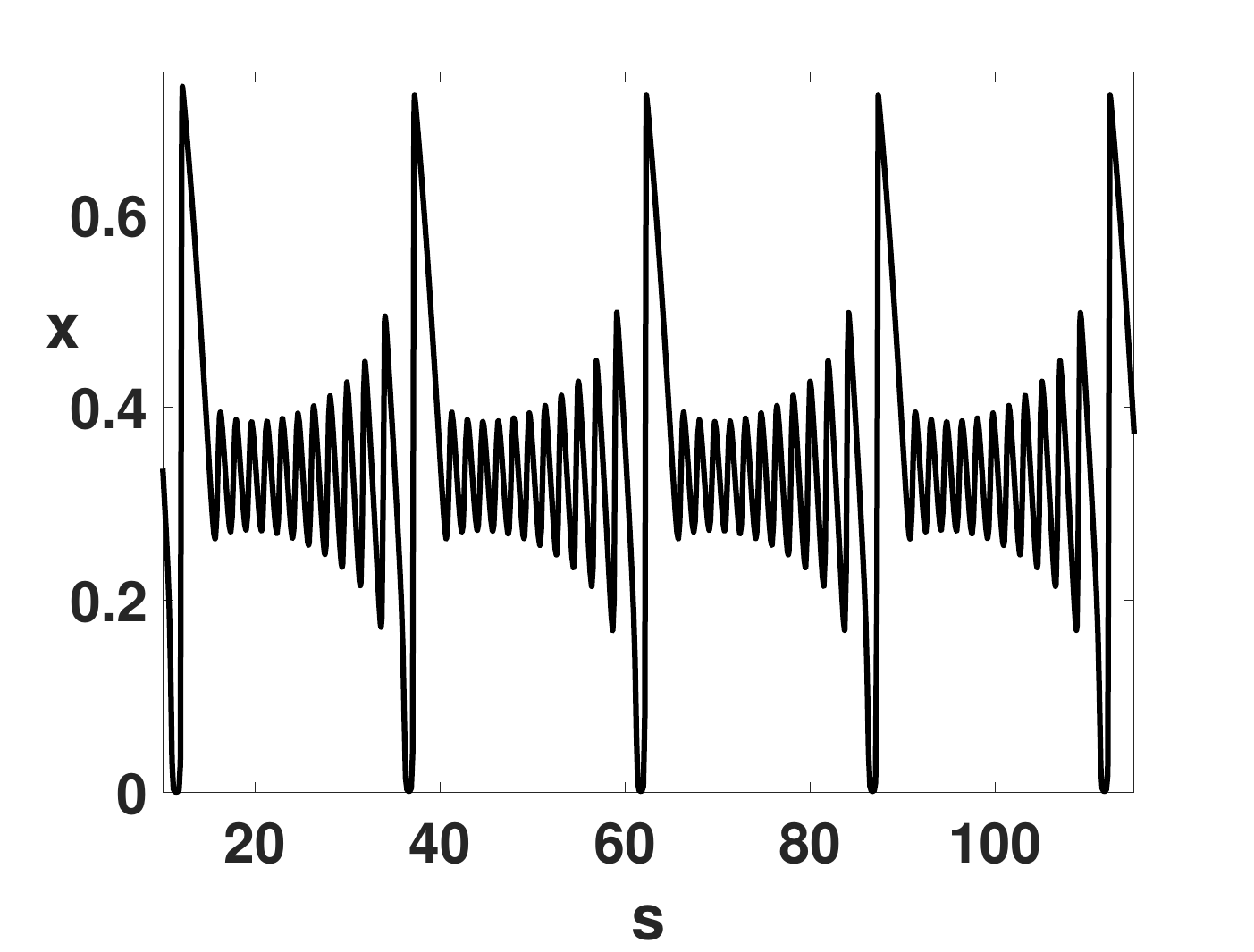}}  
  \caption{(A)-(B): Phase space and time series of an MMO orbit of signature $1^{11}$ for $h=0.633$. The SAOs rotate around the weak eigendirection (in red) and are induced by the rotational properties of the weak canard (in blue). Also shown is the folded node singularity (in cyan).  Here $\beta_1=0.1977$ and the other parameters are as in (\ref{parvalues}).}
  \label{orbits_near_FH_MMO}
\end{figure}

System (\ref{nondim3}) exhibits MMOs for different combinations of parameter values. Figure \ref{orbits_near_FH_MMO} represents an MMO orbit corresponding to  $\beta_1=0.1977$ and $h=0.633$ and other parameter values as in (\ref{parvalues}). These parameter values are chosen from a regime near a co-dimension 2 bifurcation point (see figure \ref{two_par}).   Note that the SAOs associated with the MMO orbit in  figure \ref{timeseries_08_819} have much smaller amplitudes than the ones seen in figure \ref{orbits_near_FH_MMO}. The small oscillations start with low amplitude, decrease and then increase before entering into the relaxation regime. In contrast, the small oscillations of the MMO orbit in figure \ref{orbits_near_FH_MMO}  start with moderate amplitude and then grow slowly with increasing amplitude. MMOs with former characteristics have been attributed to canards at folded nodes, while the latter characteristics have been attributed to singular Hopf points in \cite{G}. However, the presence of a  folded node near the singular Hopf point in this system (see figure \ref{desing}) may allow MMOs with only increasing-in-amplitude oscillations but not decreasing-in-amplitude, and hence we see both types of features in this model.
 Similar characteristics have been observed in other models in the literature such as in \cite{RWK}.

\subsection{MMOs as prolonged transients}

Considering $\beta_1=0.25$ and other parameter values as in (\ref{parvalues}), we noted that a supercritical Hopf bifurcation of the coexistence equilibrium occurs at $h\approx 0.7803$, giving birth to a stable limit cycle $\Gamma_h$. 
\begin{figure}[h!]     
  \centering 
{\includegraphics[width=14.0cm]{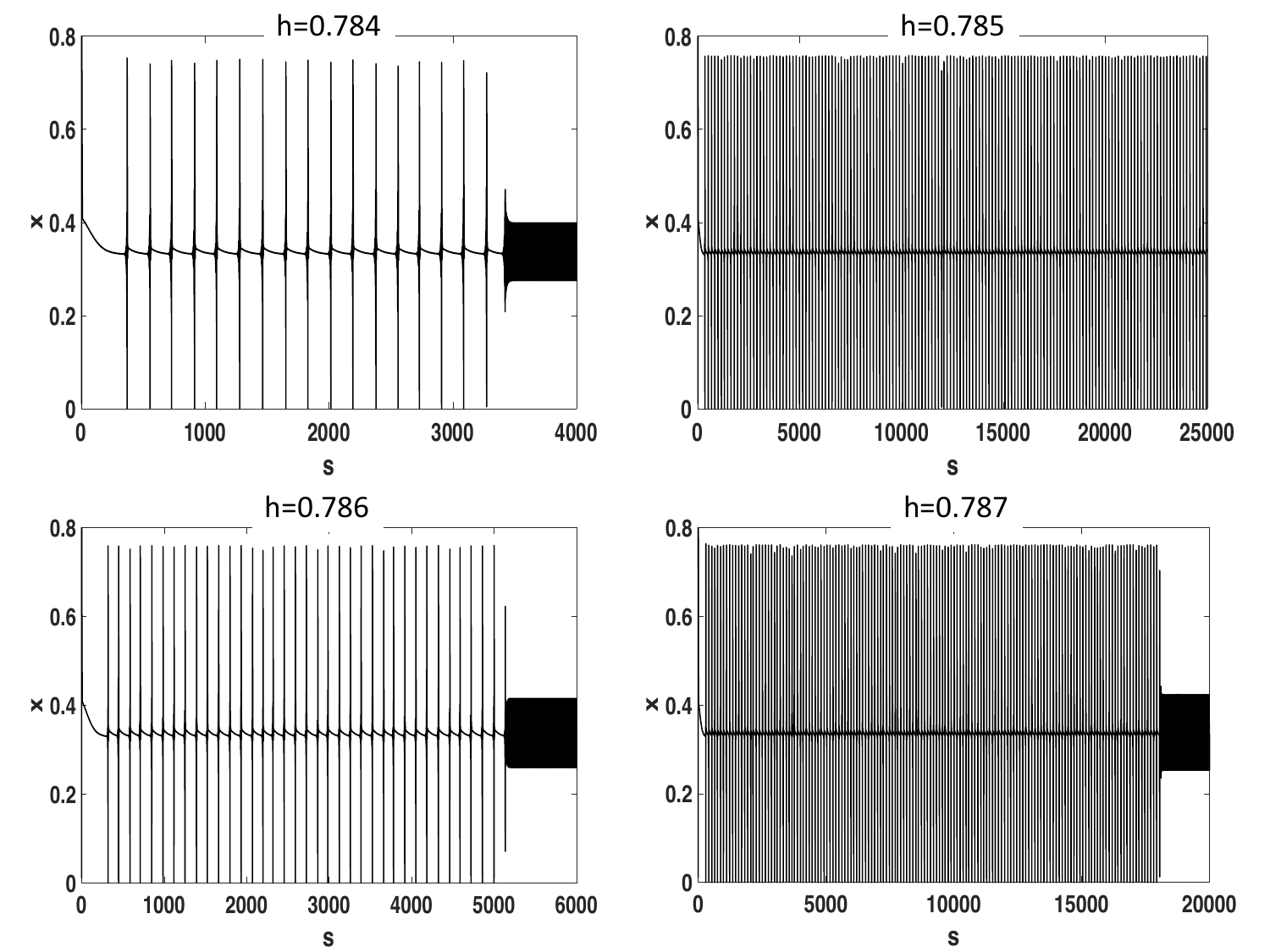}}
 \caption{Time series in $x$ demonstrating long transient dynamics in form of MMOs  past supercritical Hopf bifurcation as $h$ is varied. The initial condition chosen is $(0.01, 0.01, 0.12)$. Note the duration of the transients varies with $h$.}
 \label{transient}
\end{figure}
Interestingly, in the parameter regime where $\Gamma_h$ is stable, i.e. for $0.7803<h<0.797$, long term transients in form of chaotic MMOs are observed before the system approaches its asymptotic state. 
The duration of the transient depends sensitively on the initial values of the state variables, and may last for a significantly long amount of time for certain initial conditions. Also, with the same initial condition,  chaotic transients with a correspondingly distinct duration are observed as the control parameter is varied as shown in figure \ref{transient}.  
\begin{figure}[h!]     
  \centering 
\subfloat[]{\includegraphics[width=7.5cm]{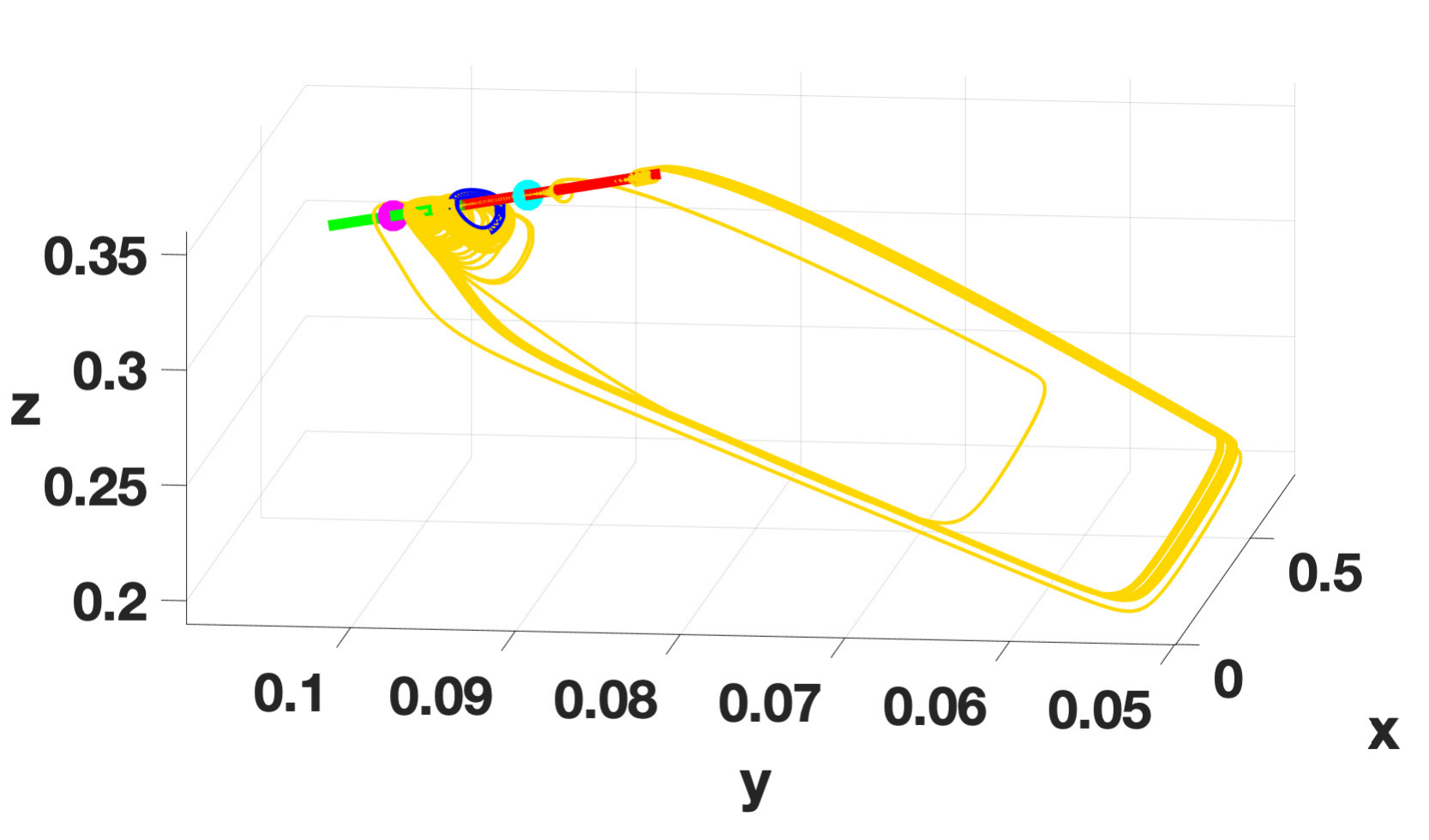}}
\qquad
\subfloat[]{\includegraphics[width=5.0cm]{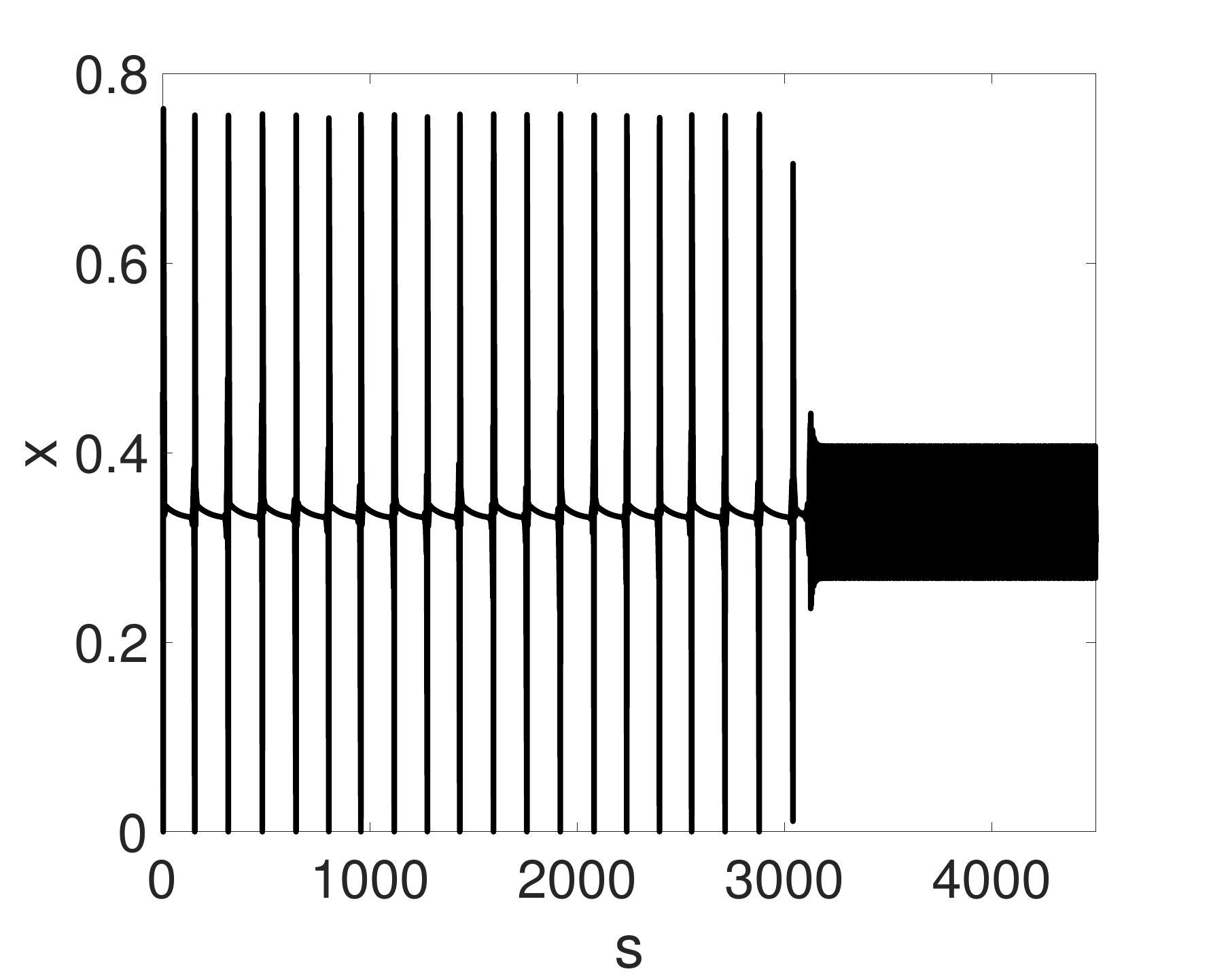}}\qquad
\subfloat[]{\includegraphics[width=7.5cm]{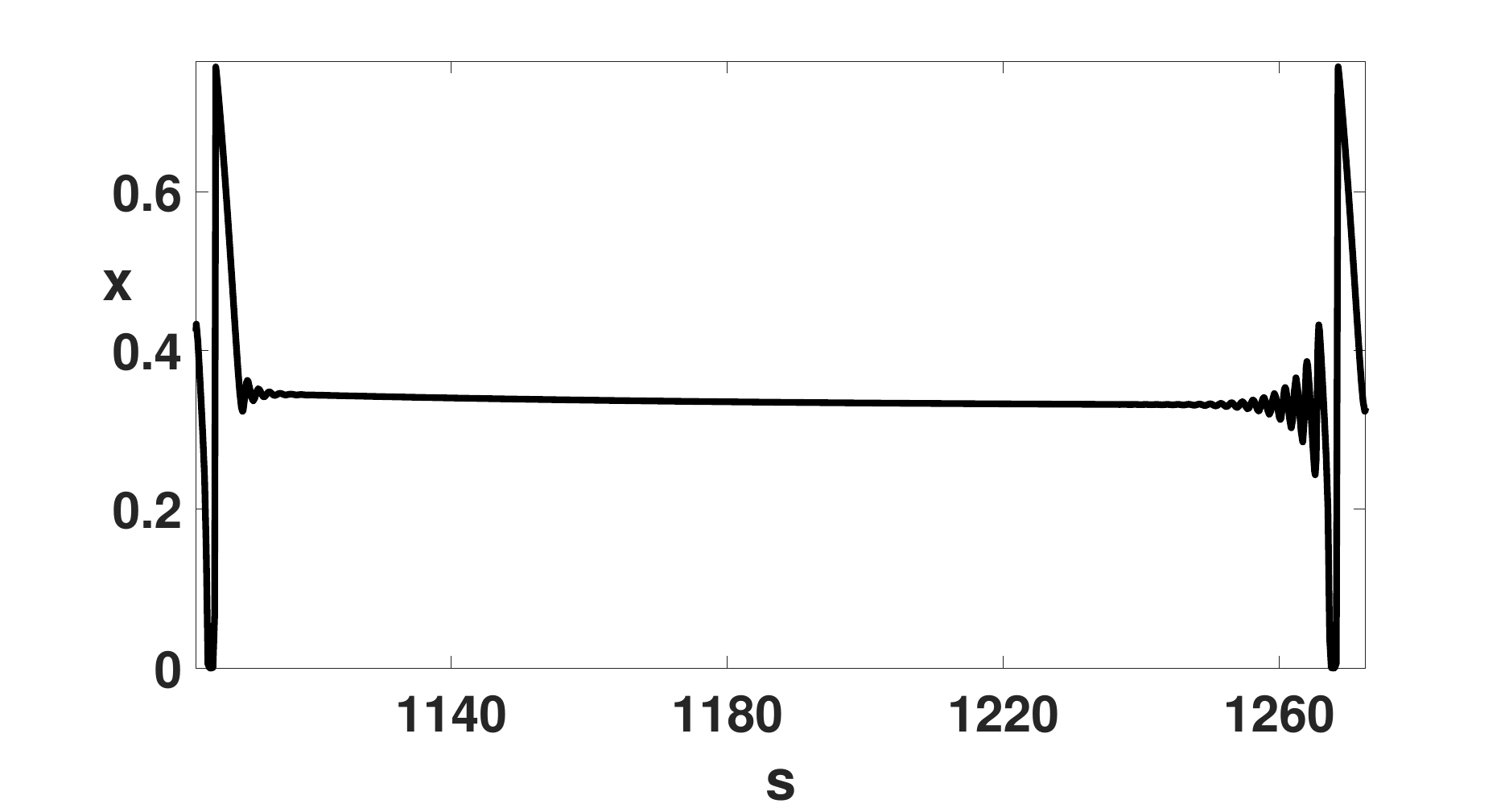}}\qquad
\subfloat[]{\includegraphics[width=7.65cm]{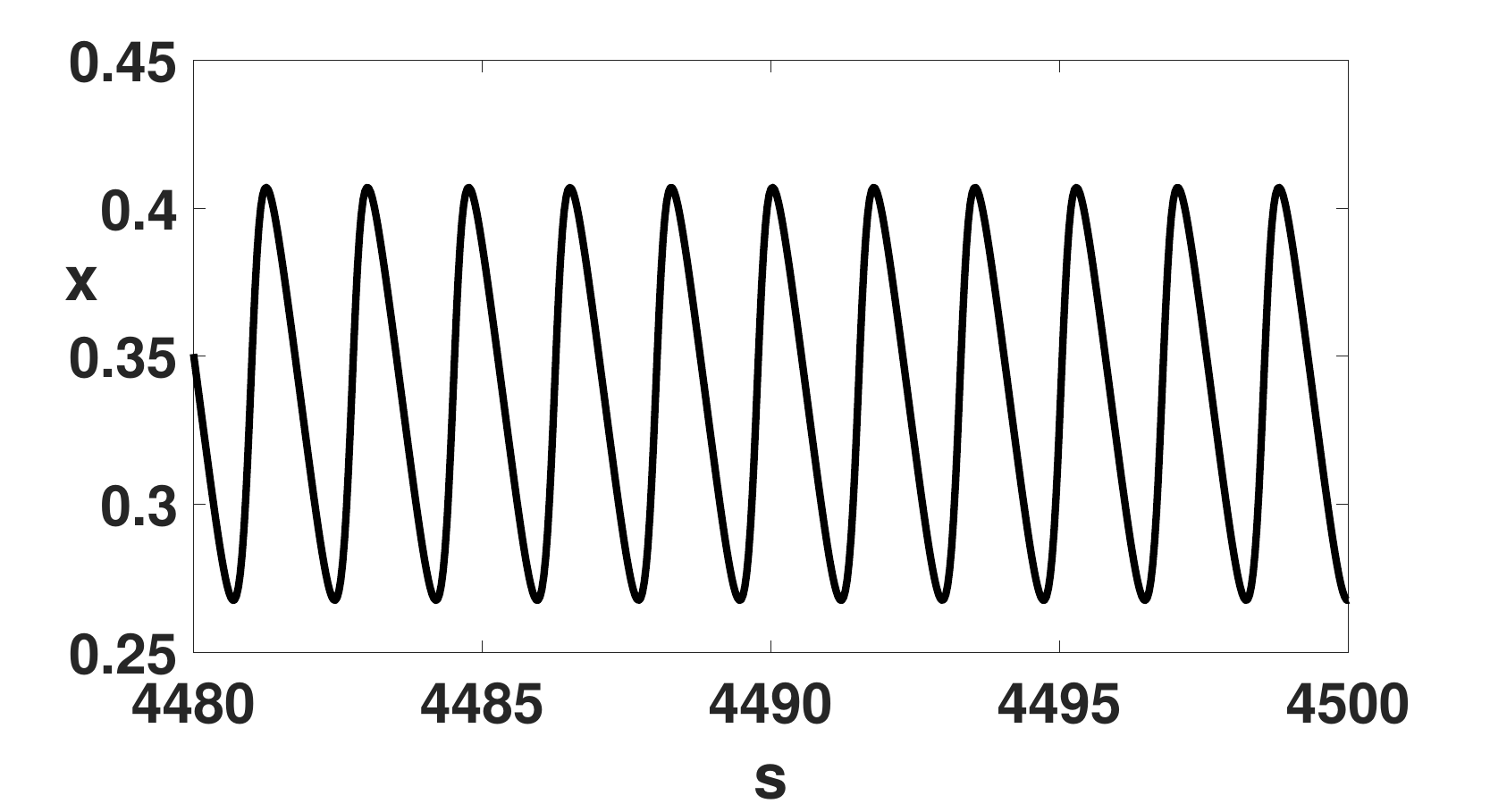}}
 \caption{(A) A trajectory exhibiting MMOs while approaching $\Gamma_{h}$ (blue).  The trajectory tightly wraps around the weak eigendirection $v_w$  (in red) and approaches $E_1^*$ (marked as a magenta dot), inducing SAOs in the MMOs. The folded node singularity is marked as the cyan dot and the stable eigendirection of $E_1^*$ in green .  (B) Corresponding time series in $x$. (C) A zoomed view of the SAOs between two LAOs. Note the amount of time spent near the folded node and the equilibrium. (D) Time series of $\Gamma_{h}$, the periodic stable attractor. Here, $\beta_1=0.25$, $h=0.785$ and  the other parameter values are as in (\ref{parvalues}).}
 \label{bistable_attract}
\end{figure}

\begin{figure}[h!]     
  \centering 
\subfloat[]{\includegraphics[width=7.7cm]{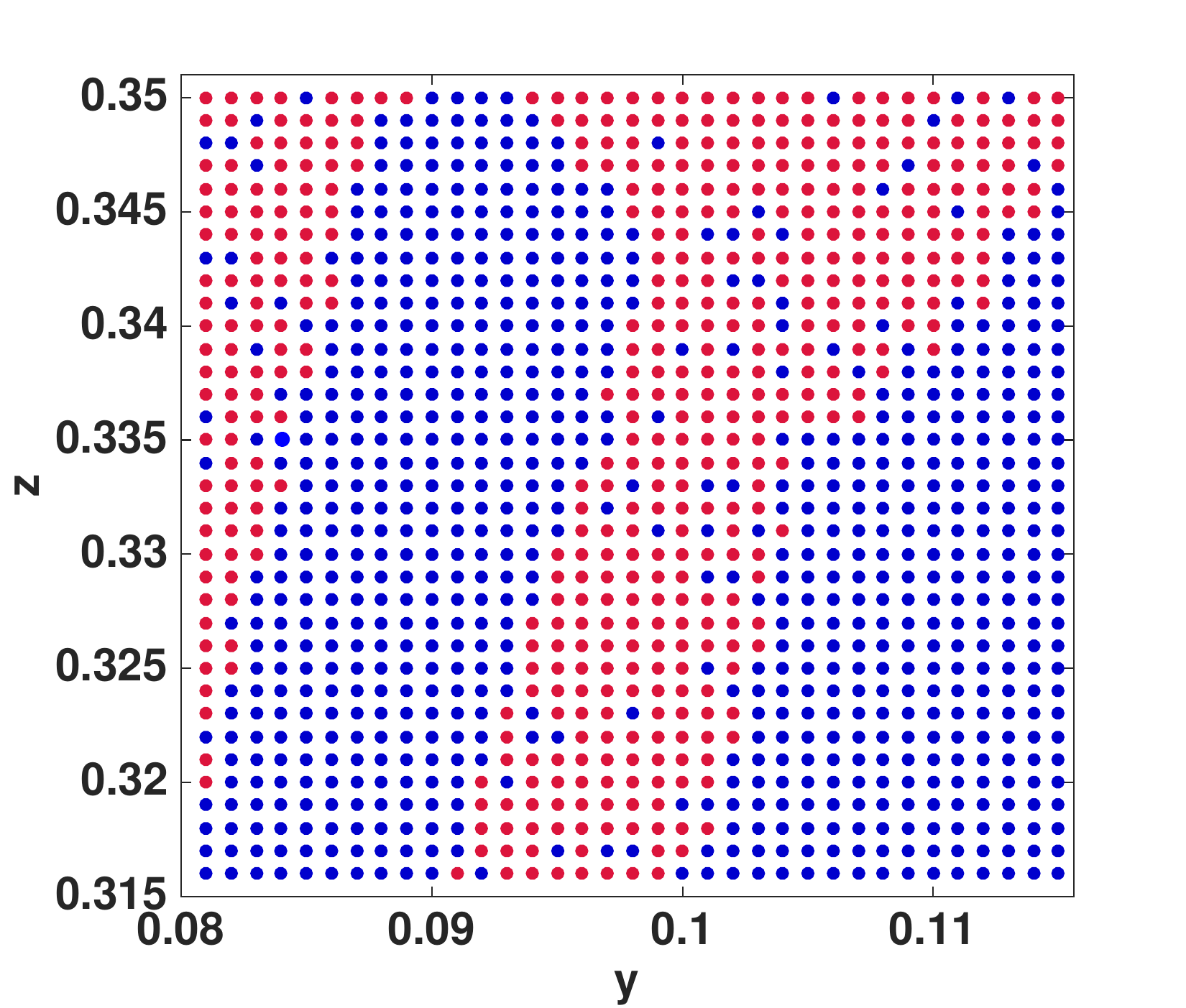}}
\quad
\subfloat[]{\includegraphics[width=7.7cm]{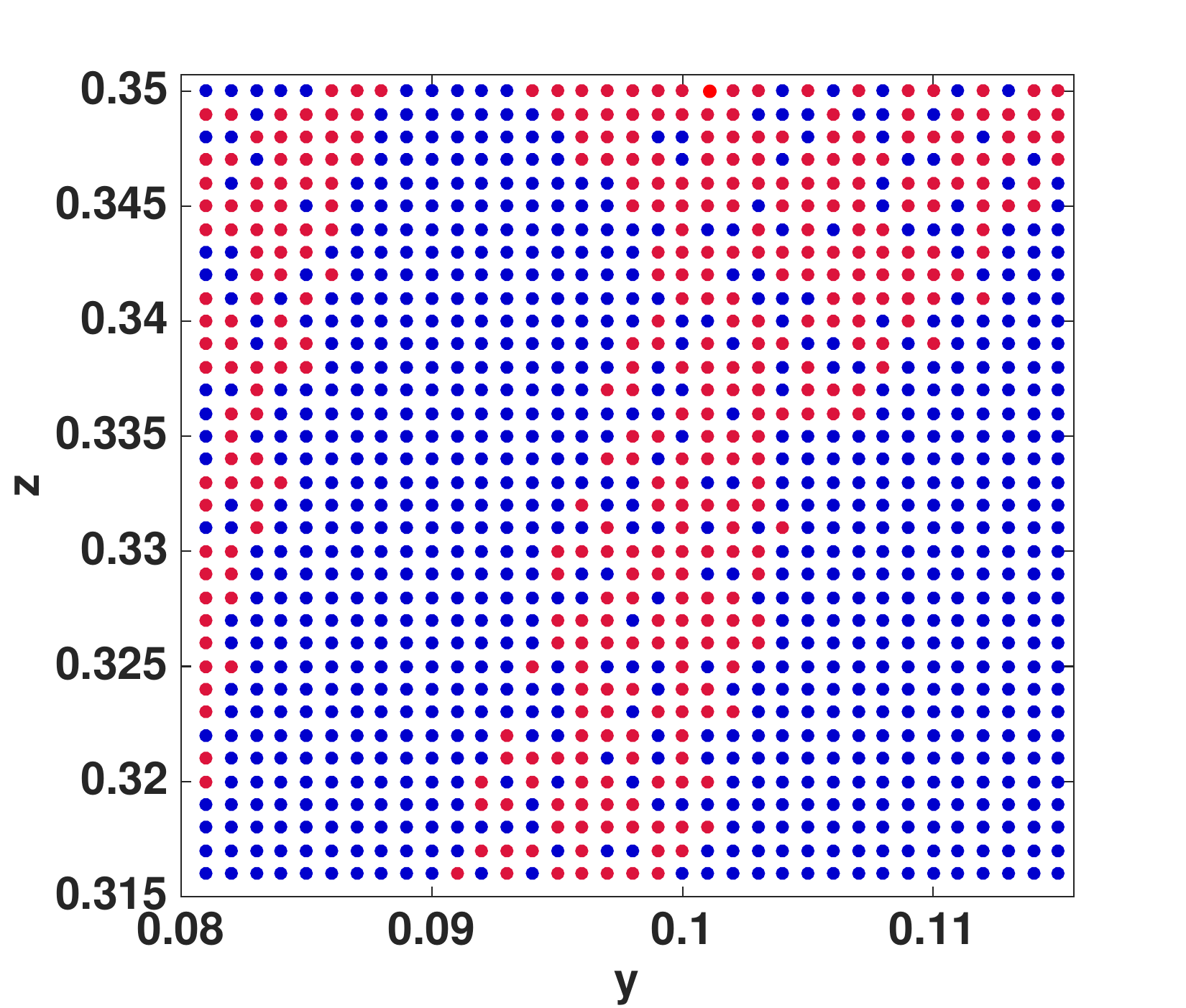}}
 \caption{Transient basins of attraction  for $h=0.785$ restricted to the plane  $x=0.3428$ computed over the interval $[0, N]$. Blue dots lie in the basin of short transients and the red dots in the basin of long transients. All the other parameter values are as in (\ref{parvalues}) with $\beta_1=0.25$.  (A) $N=2000$ (B) $N=5000$.}
 \label{trans_mmo_sao}
\end{figure}

 Figure \ref{bistable_attract} demonstrates the transition from chaotic MMOs to the asymptotic dynamics, i.e. the Hopf limit cycle, for a fixed set of parameter values.  Such complex oscillatory patterns in population densities persist for thousands of generation  before the system undergoes a sudden change to  a completely different state. A ``transient basin" can be thus defined as the set of initial conditions which produce chaotic transients whose duration is within a given interval $[0, N]$ \cite{SCVB}. Numerical approximations of transient basins for a fixed value of $h$ are plotted in figure \ref{trans_mmo_sao} by choosing equally spaced points on a  grid of size $[0.08 , 0.115] \times [0.325, 0.36]$ on the plane $x=0.3428$ and integrating trajectories starting at each point on this grid. Points on the grid which lead to transients of duration less than $N$ are plotted in blue, and are  referred to as ``basins of short transients". On the other hand, points in red will be referred to as ``basins of long transients". The transient basins are densely intermixed as can be seen in figure \ref{trans_mmo_sao}.

\subsection{Relaxation oscillations as long term transients}
Another interesting feature observed in system (\ref{nondim3}) is the existence of relaxation oscillations as prolonged transients before the system approaches the boundary equilibrium state $E_{xz}$. This occurs when one of the predators, namely $z$ is assumed to be more efficient than $y$. To this end, we let $\beta_1=0.1923$ and consider parameter values as in (\ref{parvalues1}). 
 \begin{figure}[h!]     
  \centering 
\subfloat[]{\includegraphics[width=7.0cm]{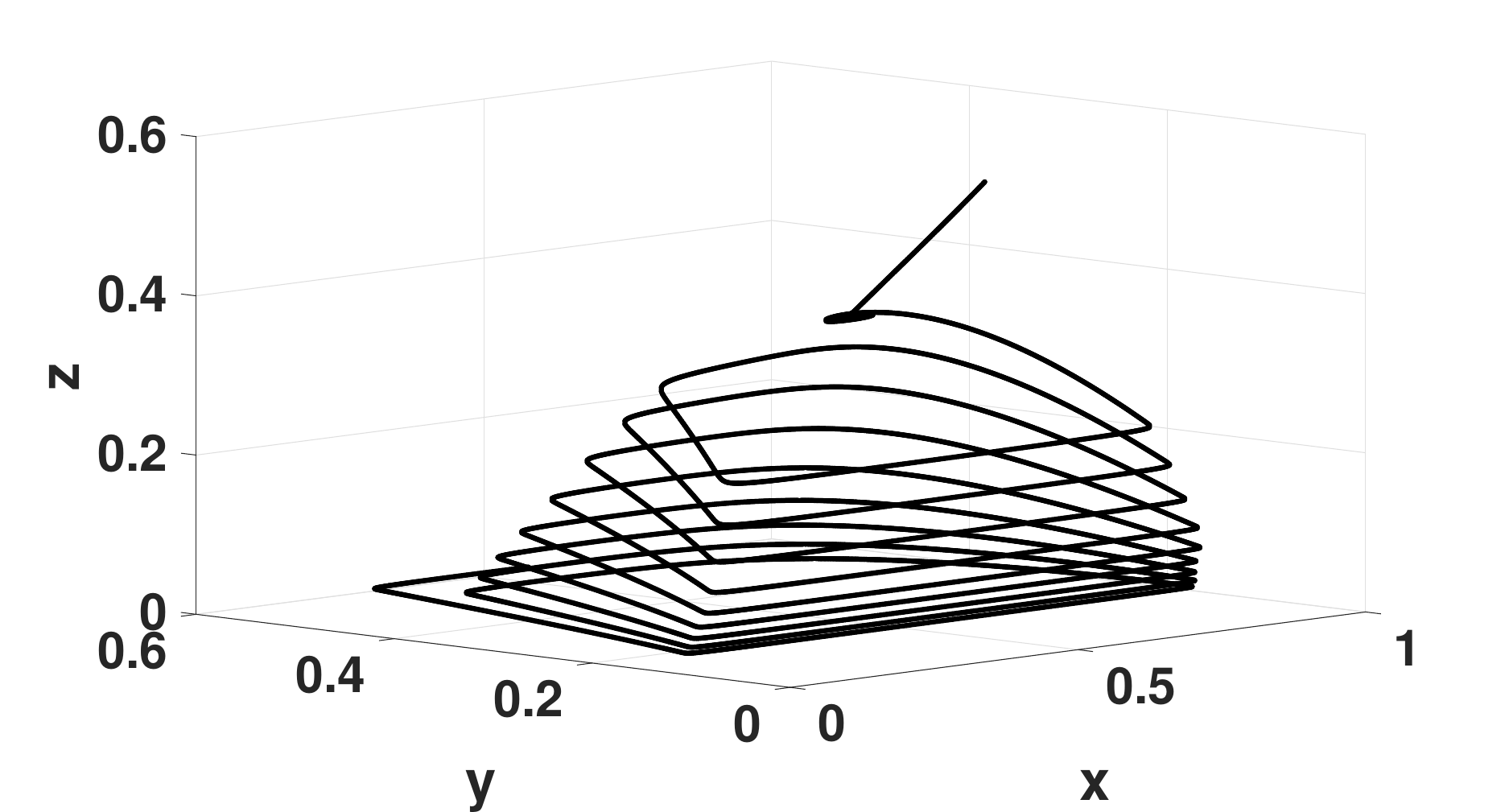}}
\quad
\subfloat[]{\includegraphics[width=7.0cm]{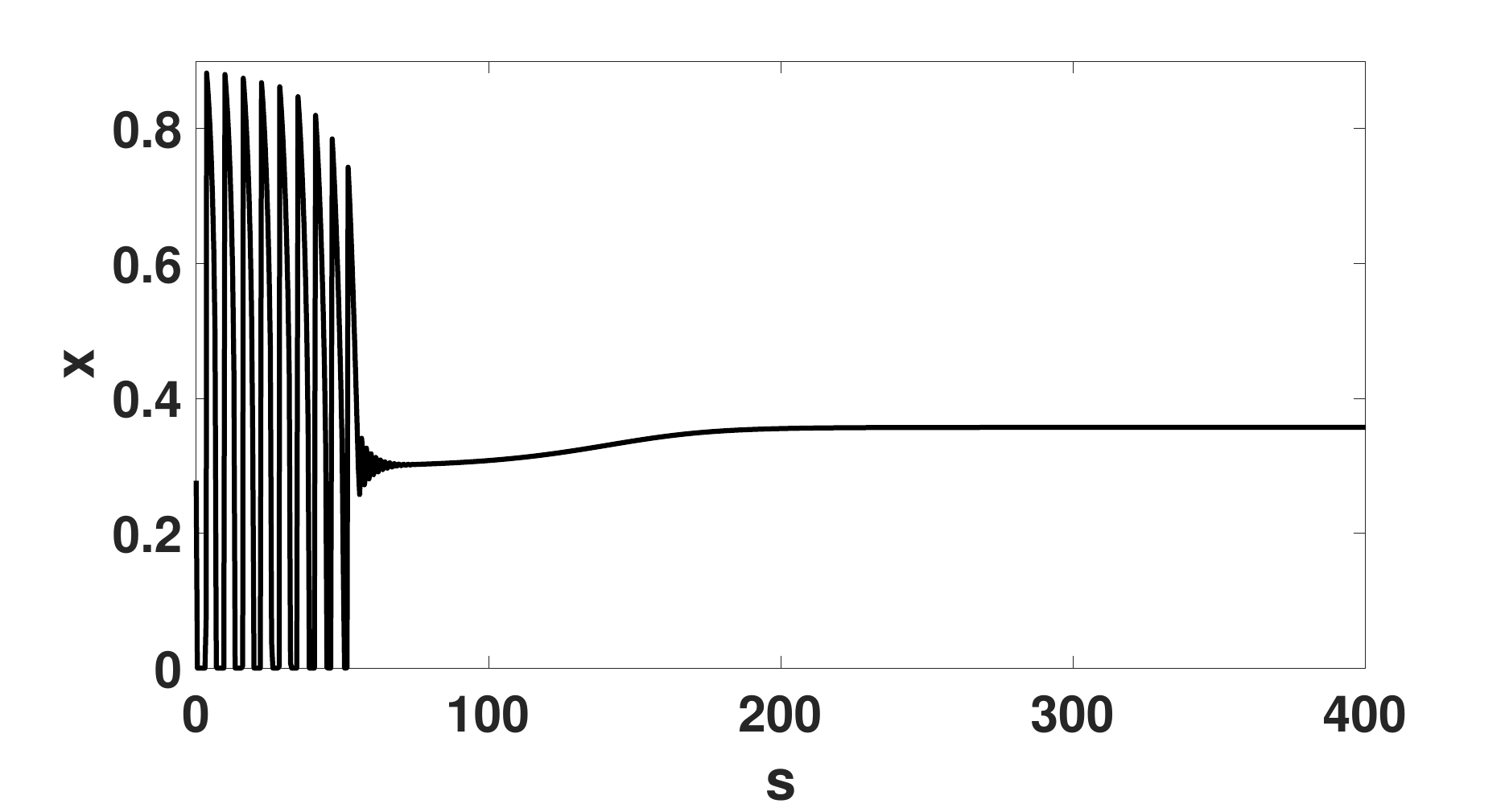}}
 \caption{Phase portrait and time series of a trajectory starting out at $(0.27795, 0.4252, 0.0598)$ exhibits relaxation oscillations as transient dynamics and eventually asymptotes to $E_{xz}=(0.357, 0, 0.615)$ for $\beta_1=0.1923$, $h=0.2649$ and other parameter values as in (\ref{parvalues1}).}
 \label{bistability}
 \end{figure}
A subcritical Hopf bifurcation of the coexistence equilibrium occurs at $h=0.2646$, followed by a supercritical torus bifurcation at $h=0.2648$. For $h>0.2648$, the semi-trivial equilibrium state $E_{xz}$ is locally asymptotically stable and coexists with the periodic orbit born out of torus bifurcation. The periodic orbit has a relatively smaller basin of attraction and the system can either approach the periodic orbit or $E_{xz}$. Transient dynamics in  form of relaxation oscillations are observed as the system approaches $E_{xz}$ as shown in figure \ref{bistability}. These types of transient dynamics which may apparently seem sustainable  on ecologically relevant timescales could result into species extinction, a property that has been studied in many ecosystems \cite{Hastings}.

\section{Bifurcation analysis}



In this section, we explore the bifurcation structure of system (\ref{nondim3})  and find other interesting dynamics that are generated by this generic model.  We begin by treating $h$ as the primary bifurcation parameter and $\beta_1$, the predation efficiency  of $y$, as the secondary parameter, keeping the other parameter values fixed at (\ref{parvalues}) or (\ref{parvalues1}). Besides studying the effect of $h$ and $\beta_1$ on the system, we will also explore the role of $\beta_2$. The two different values of $\beta_2$  considered in (\ref{parvalues}) and (\ref{parvalues1}) will elucidate the influence of predation efficiency of $z$ on the system. We will compute two-parameter bifurcation diagrams in $(h, \beta_1)$ and then compare the bifurcation structures of system (\ref{nondim3})  for the two sets of parameters to see  how the structure changes with the predation efficiency of $z$.

\subsection{A two-parameter bifurcation of system (\ref{nondim3}) with similar predation efficiencies.} 
 \begin{figure}[h!]     
  \centering 
  {\includegraphics[width=10.5cm]{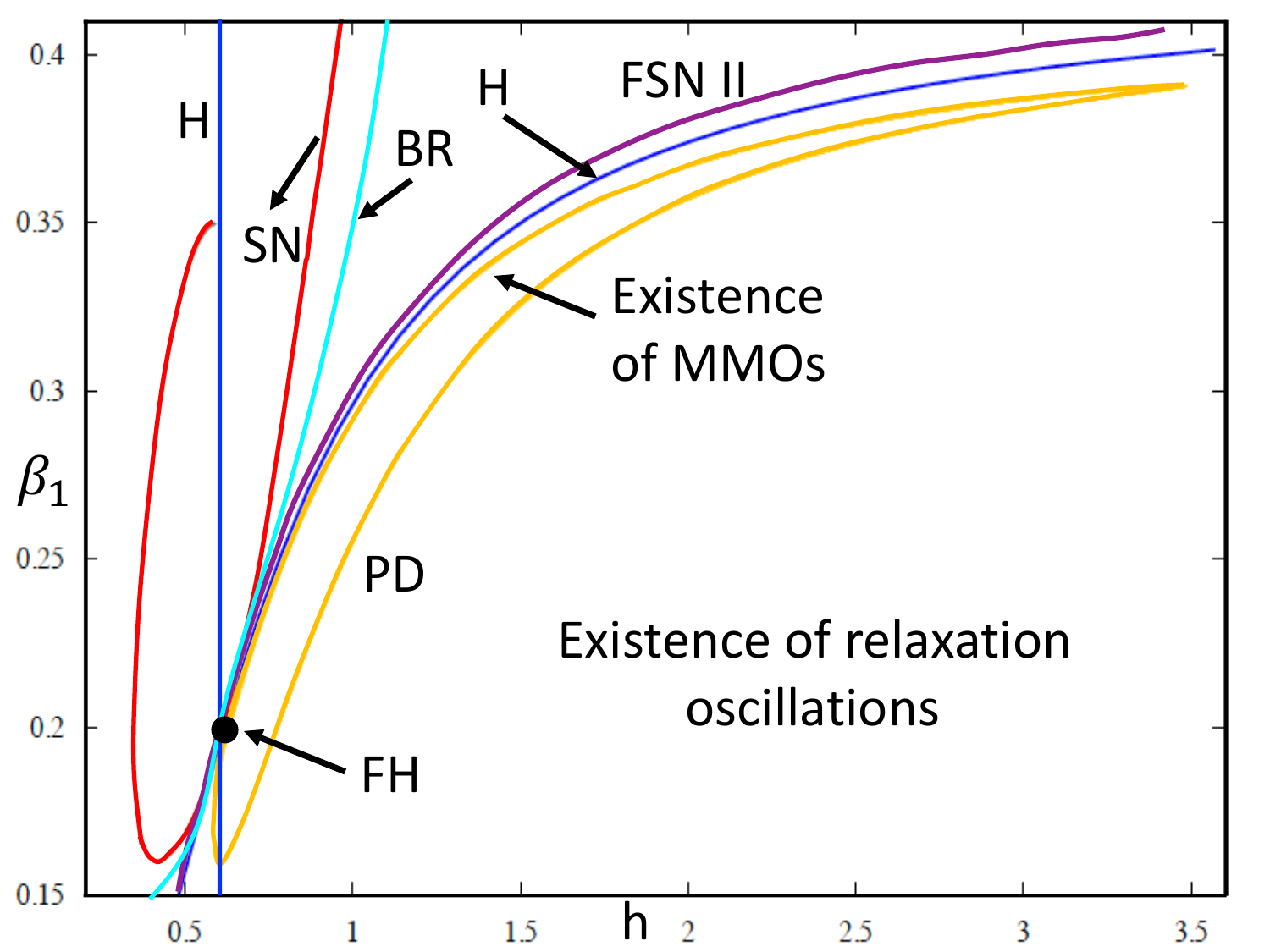}}
  \caption{A two-parameter bifurcation diagram in $(h, \beta_1)$ space. The other parameter values are as in (\ref{parvalues}). FH: fold-Hopf bifurcation (marked by black dot),  SN: saddle-node bifurcation (red curve), H: Hopf bifurcation (blue curve), FSN II bifurcation (purple curve), PD: period-doubling bifurcation (yellow curve), BR: branch curve (cyan). The equilibria lying to the right of BR are biologically feasible. The plot was generated in XPPAUT \cite{E}.}
  \label{two_par}
\end{figure}

Keeping all the parameter values fixed as in (\ref{parvalues}), we consider a two-parameter bifurcation in $(h, \beta_1)$ parameter space. The bifurcation diagram is shown in figure \ref{two_par}. The  equilibrium point $E_{xz}$, lying on the $xz$-plane, undergoes a supercritical Hopf bifurcation along the vertical Hopf line $H$ (shown in blue), and the coexistence  (non-trivial) equilibrium undergoes a supercritical Hopf bifurcation along the Hopf curve, also denoted by $H$. The FSN II curve and the Hopf curve lie very close to each other.  A saddle-node bifurcation of the non-trivial equilibria of (\ref{nondim3}) occurs along the SN curve (shown in red). The equilibria lying on the SN curve may or may not be biologically feasible (as they could have one or more negative components). In between the SN curve and the branch curve (BR), system (\ref{nondim3}) could have exactly two or no nontrivial equilibria.   A unique nontrivial equilibrium of (\ref{nondim3}) exists to the right of BR. In the region lying in between BR and the Hopf curve, the nontrivial equilibrium is locally asymptotically stable. To the right of the Hopf curve, oscillatory dynamics such as MMOs and relaxation oscillations appear. The PD curve, corresponding to the period-doubling bifurcation of the periodic orbit born out of the Hopf bifurcation demarcates MMOs from relaxation oscillations in the  parameter space. The branch curve (BR)  meets with the Hopf curve at the fold-Hopf bifurcation point (FH) which occurs at $(\beta_1, h)=(0.197019, 0.60803)$.  Near the fold-Hopf bifurcation, which is a  co-dimension two bifurcation  \cite{GH, K1},  other bifurcations may occur including saddle-node bifurcations of periodic orbits and  bifurcation of Shil'nikov homoclinic orbit to a saddle-focus. The latter type of bifurcation can serve as an organizing center for the MMOs near the fold-Hopf point, and will be a subject for future investigation. 


\subsubsection{Dynamics near the fold-Hopf point} 
\begin{figure}[h!]     
  \centering   
      \subfloat[]{\includegraphics[width=7.5cm]{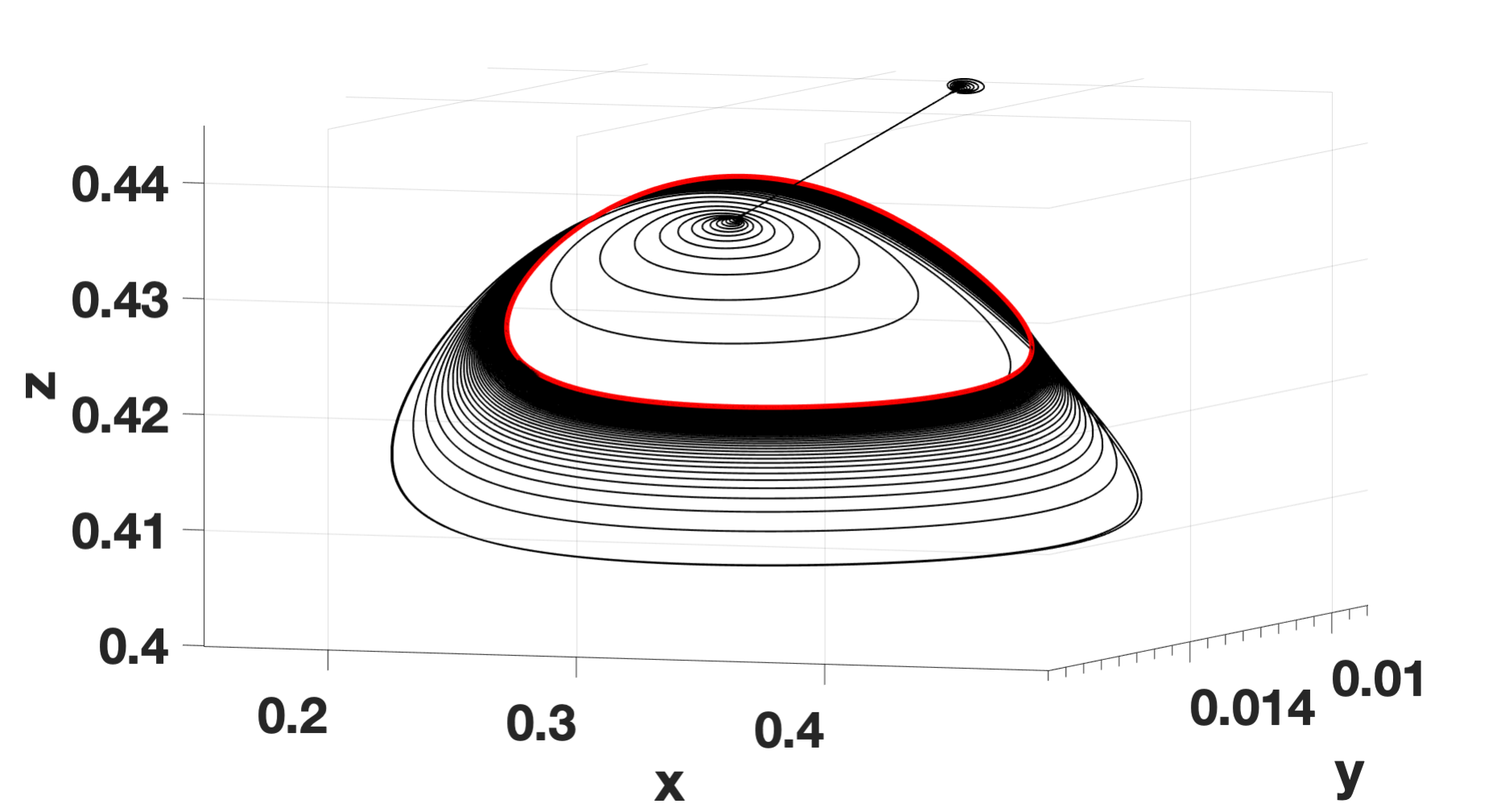}}\qquad
         \subfloat[]{\includegraphics[width=5.7cm]{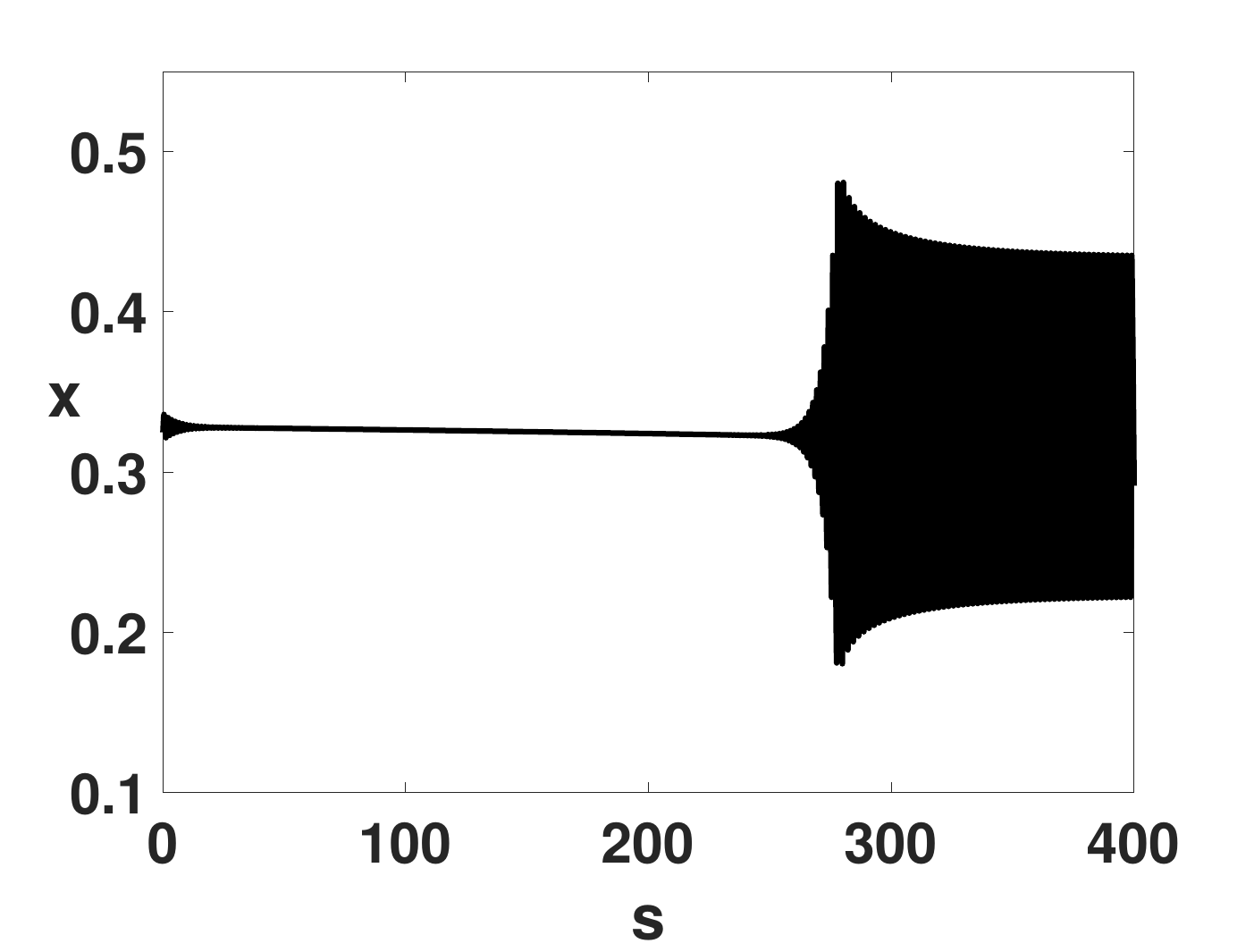}}\qquad
  \caption{(A)-(B): Phase space and time series of a trajectory for $h=0.615$. The trajectory approaches  $E_1^*$ along $W^s(E_1^*)$ and spirals out along $W^u(E_1^*)$ as it eventually settles down to the periodic orbit (red), which emerges through a Hopf bifurcation.  Here $\beta_1=0.1977$ and the other parameters are as in (\ref{parvalues}).
   }
  \label{orbits_near_FH}
\end{figure}
In a vicinity of the FH point, above the Hopf curve $H$ and below the branch curve BR in figure \ref{two_par}, system (\ref{nondim3}) possesses two coexistence equilibria $E_1^{*}$ and $E_2^*$. The equilibrium $E_1^{*}$ is a stable focus-node and $E_2^{*}$ is a saddle-focus with two-dimensional stable manifold and one-dimensional unstable manifold. As the control parameter $h$ is increased, $E_1^{*}$ undergoes a supercritical Hopf bifurcation and  thereafter becomes a saddle-focus with one-dimensional stable manifold, $W^s(E_1^*)$, and two-dimensional unstable manifold, $W^u(E_1^*)$. The equilibrium $E_2^*$ no longer lies in the first octant and the system asymptotes to a periodic orbit that emerges from the Hopf bifurcation of $E_1^*$ (see figure \ref{orbits_near_FH}). On further increasing $h$, MMOs (both aperiodic and periodic orbits) appear (see figure \ref{orbits_near_FH_MMO}). The SAOs associated with the MMOs are {\emph{canard induced}}. Such dynamics will be discussed in the next section. 
 The amplitude of the SAOs depend sensitively on $\beta_1$, which will be also explored in this paper. 

\subsubsection{Horizontal bifurcations} Treating the two predators to be equally efficient, i.e. by letting  $\beta_1=0.35$  and the other parameter values as in  (\ref{parvalues}), a similar bifurcation diagram as in figure \ref{one_par_xpp} is obtained. A supercritical Hopf bifurcation of the coexistence equilibrium occurs at $h \approx 1.5$. MMO orbits with signature $1^s$, where $s\leq 5$ are more commonly seen. An MMO orbit of signature $1^5$ is shown in figure \ref{timeseries_example}. 
 \begin{figure}[h!]     
  \centering 
  {\includegraphics[width=7.5cm]{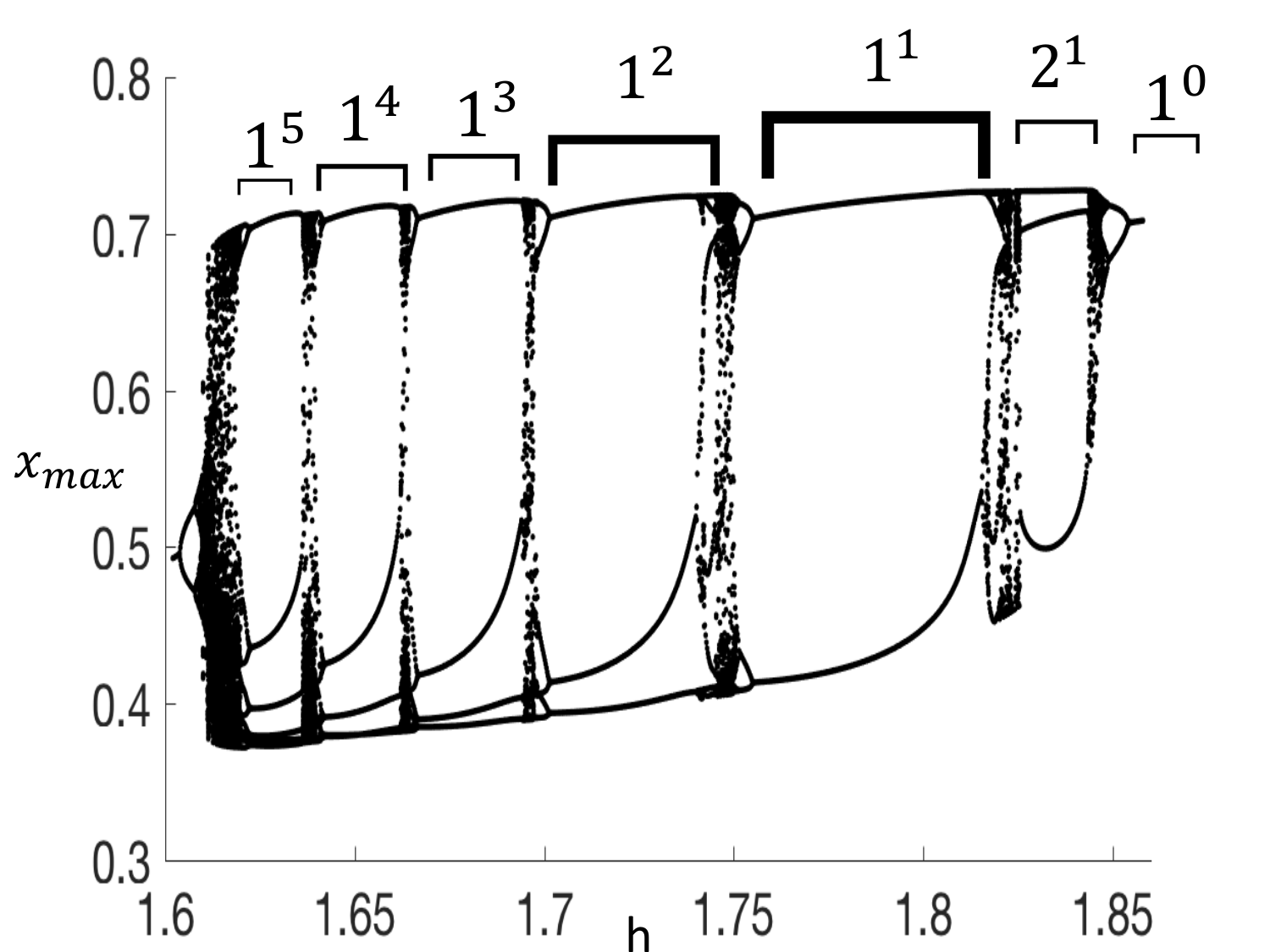}}
  \caption{A one-parameter bifurcation diagram in $h$ with $\beta_1=0.35$ and the other parameter values are as in (\ref{parvalues}). The $y$-axis represents the maximum value of $x$.}
  \label{one-par-bif-varying-beta1}
\end{figure}
Figure \ref{one-par-bif-varying-beta1} illustrates the sequence of MMOs that occurs in this case. The bifurcation diagram demonstrates how the number of small oscillations between two large oscillations vary with the intraspecific competition. Since the return time between successive population outbreaks (to a large extent) can be determined by the number of the SAOs in an MMO orbit, and therefore by the epoch of the SAOs, we note that the frequency of extreme events of outbreaks increases with the intraspecific competition amongst the predators. This reflects an important biological feature.

\subsection{A two-parameter bifurcation of system (\ref{nondim3}) with $z$  considerably more efficient than $y$.} 
\begin{figure}[h!]     
  \centering 
  {\includegraphics[width=10.5cm]{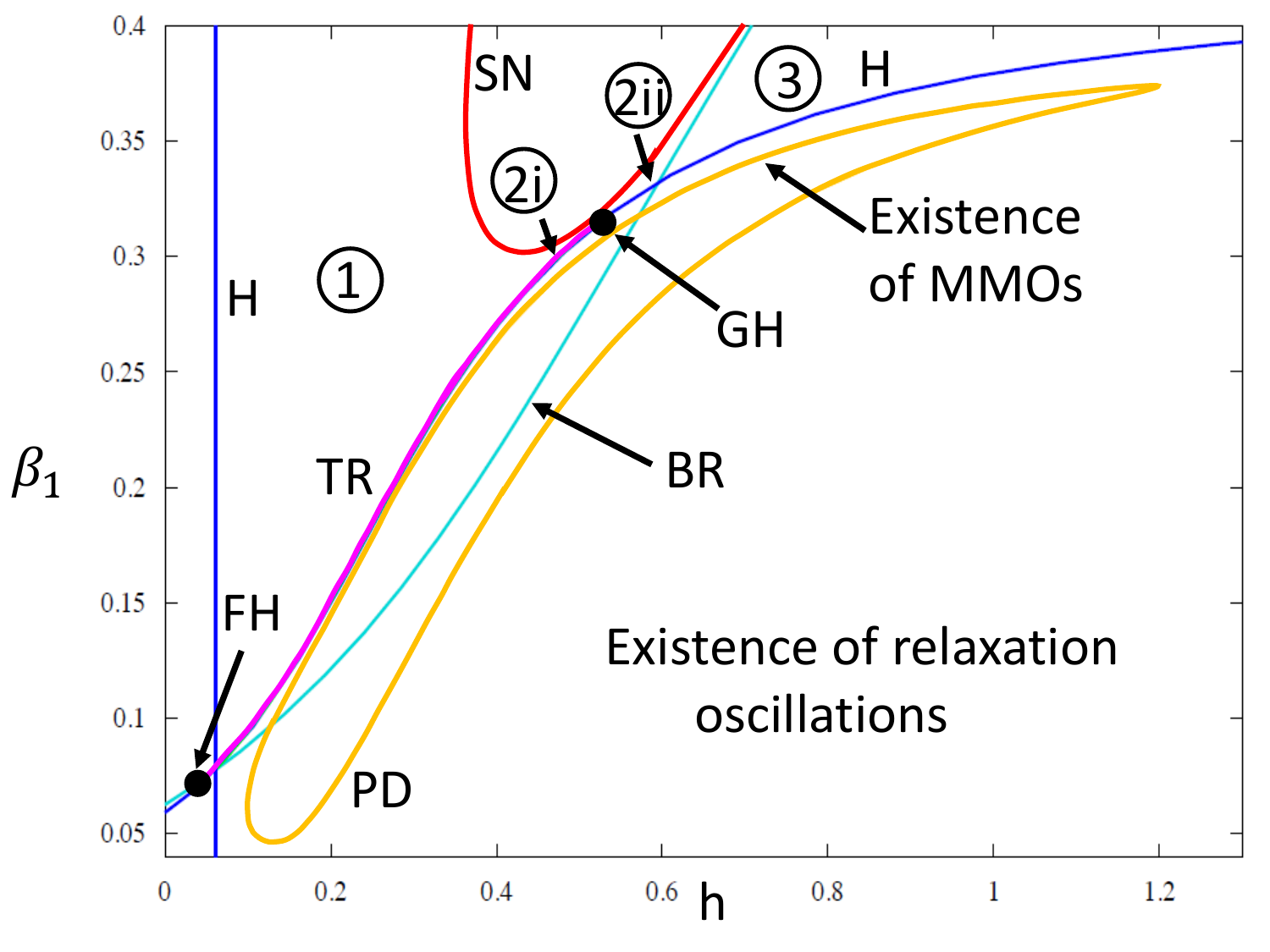}}
  \caption{A two-parameter bifurcation diagram in $(h, \beta_1)$ space. The other parameter values are as in (\ref{parvalues1}).  FH: fold-Hopf bifurcation (zero-pair),  SN: saddle-node bifurcation of equilibria (red curve), H: Hopf bifurcation (blue curve), TR: torus bifurcation (magenta curve),  BR: branch point bifurcation (cyan curve), GH: generalized Hopf bifurcation. Region 1: a unique coexistence equilibrium exists, which is unstable, Region 2i: existence of two coexistence equilibria, both unstable, Region 2ii: two coexistence equilibria exist, one of which is locally asymptotically stable, Region 3: a unique coexistence equilibrium exists, which is locally asymptotically stable. The plot was generated in XPPAUT \cite{E}.}
  \label{two_par_beta2_6}
\end{figure}

 In this subsection, we study the dynamics generated by system (\ref{nondim3}) by considering a higher value of the predation efficiency of $z$.  Keeping the other parameter values fixed as defined by (\ref{parvalues1}), we again consider a two-parameter bifurcation of system (\ref{nondim3}) in $(h, \beta_1)$ parameter space. The bifurcation diagram is shown in figure \ref{two_par_beta2_6}.  Comparing it with figure \ref{two_par}, we note that the two bifurcation diagrams are similar in many aspects, but the latter contains additional complex bifurcation points such as the generalized Hopf bifurcation. 
 
\begin{figure}[h!]     
  \centering 
  \subfloat[]{\includegraphics[width=6.75cm]{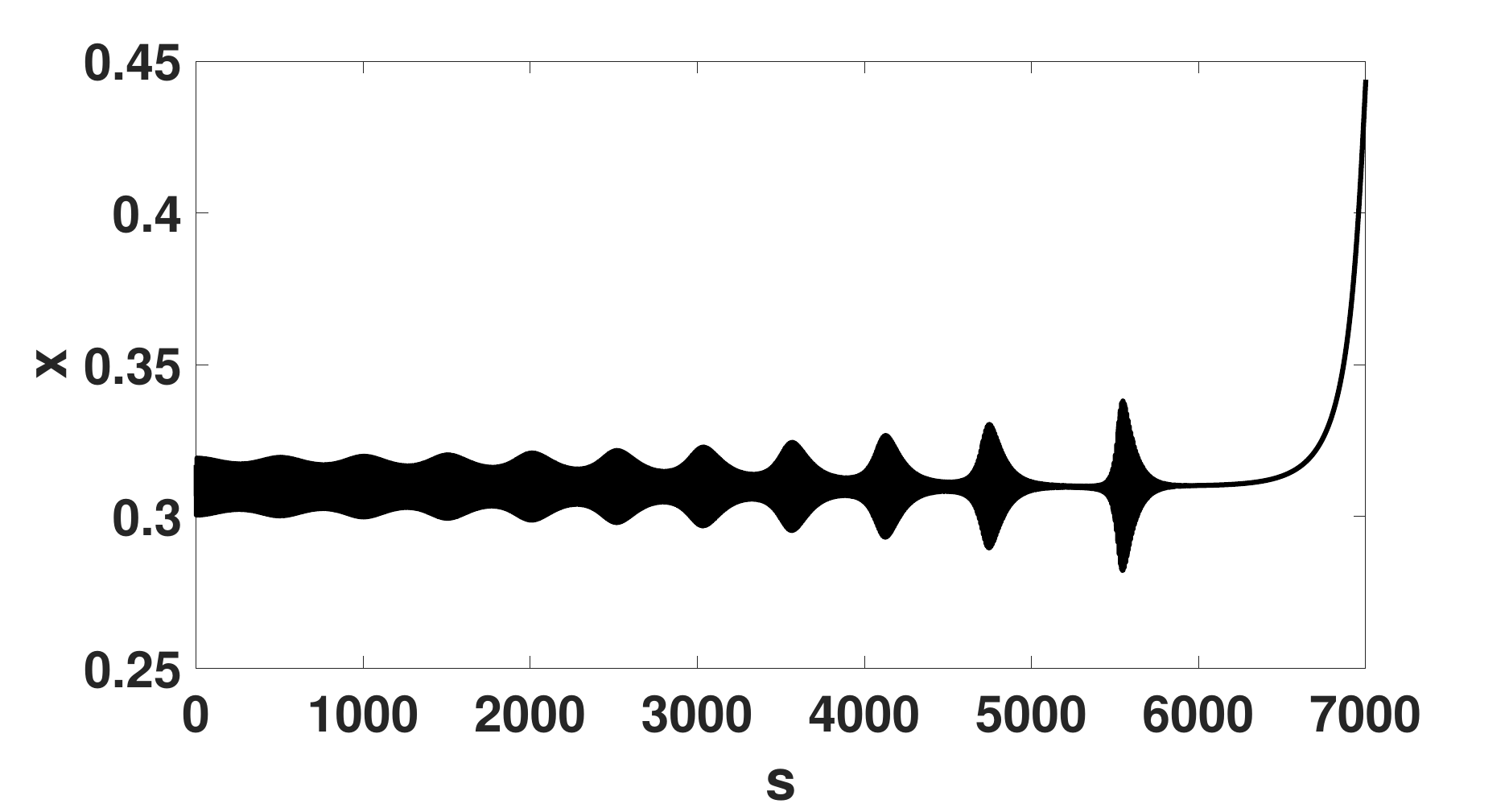}}\qquad
\subfloat[] {\includegraphics[width=6.75cm]{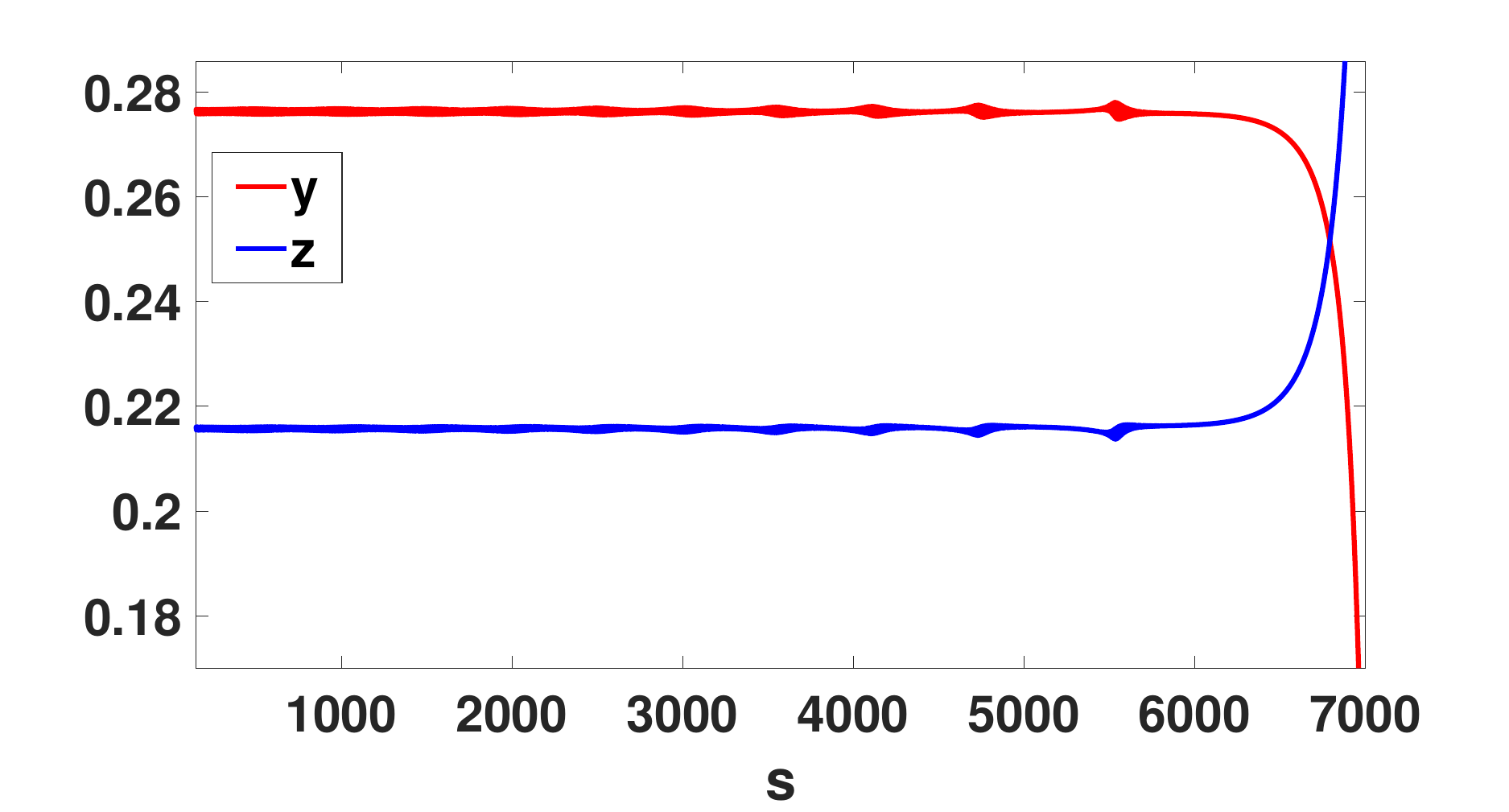}}
  \caption{Transient temporal dynamics of system (\ref{nondim3}) asymptotically approaching to the equilibrium $E_{xz}=(0.5225, 0, 0.536)$. Here $\beta_1=0.3001$, $h=0.4767$ and other parameter values as in (\ref{parvalues1}). The chosen parameters lie in region 2i in figure \ref{two_par_beta2_6}.}
  \label{timeseries_torus_bif}
\end{figure}

The equilibrium $E_{xz}$ on the $xz$ plane is locally asymptotically stable
in  regions 1 and 2i in figure \ref{two_par_beta2_6} and undergoes a supercritical Hopf bifurcation along the vertical line denoted by $H$. In a small parameter regime between the BR and SN curves, denoted by regions 2i and 2ii,  the system has two coexistence equilibria $E_1^*$ and $E_2^*$. In region 2i, the equilibrium $E_2^*$ undergoes a subcritical Hopf bifurcation as it crosses the Hopf curve $H$, while in region 2ii,  $E_1^*$ undergoes a supercritical Hopf bifurcation as it crosses the Hopf curve. A co-dimension 2 generalized Hopf bifurcation, denoted by GH,  occurs where the subcritical and the supercritical Hopf curves meet. 
In a vicinity of the subcritical Hopf curve, a torus bifurcation curve,  TR, emanates from GH and gives rise to periodic oscillations.  The TR curve terminates at the fold-Hopf (FH) bifurcation. For parameter values near the TR curve in region 2i, long lasting transient dynamics with a peculiar approach to $E_{xz}$ are seen. The temporal development of the species are shown in figure \ref{timeseries_torus_bif}. 
 \begin{figure}[h!]     
  \centering 
  \subfloat[]{\includegraphics[width=6.75cm]{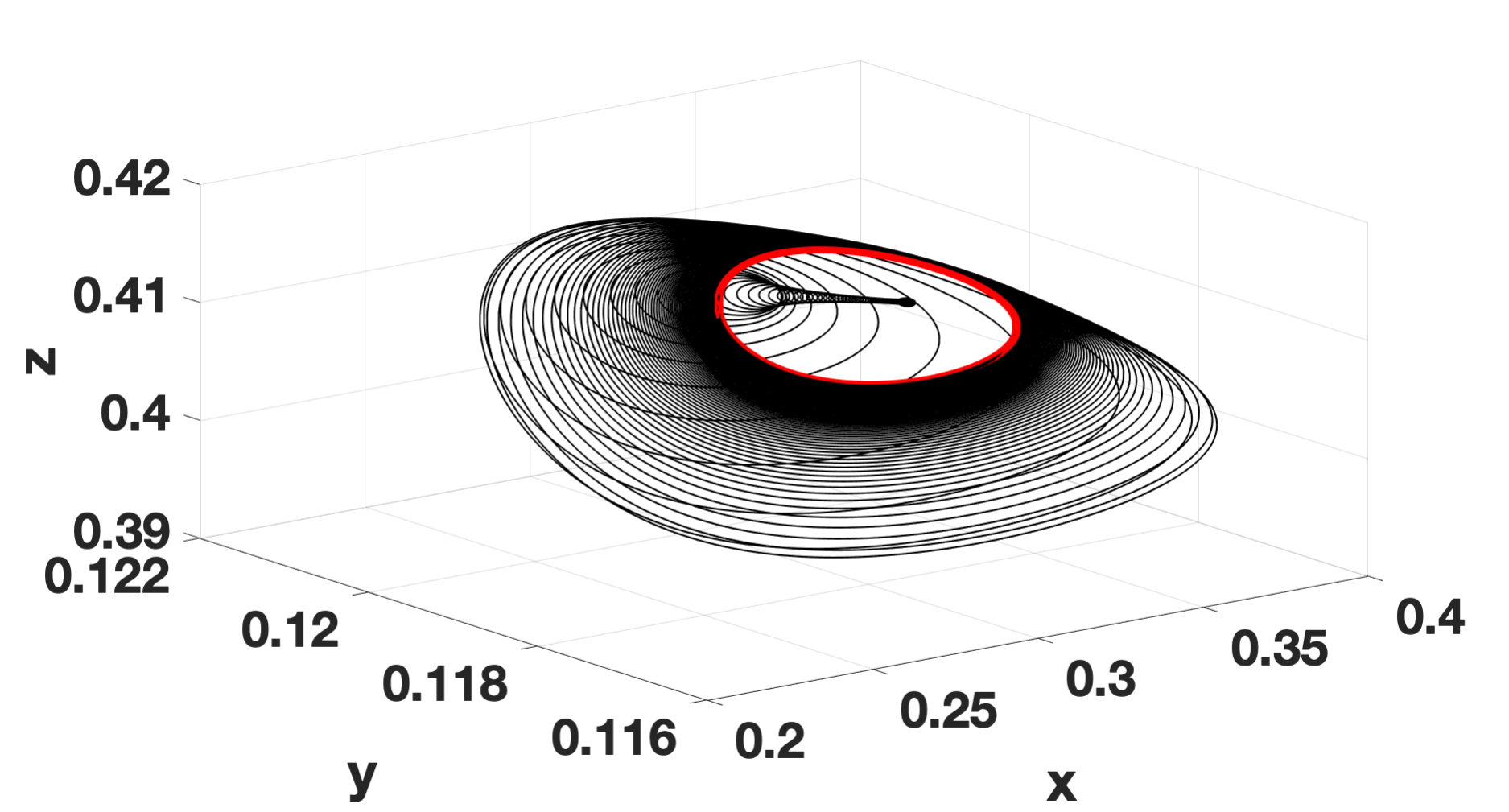}}\qquad
\subfloat[] {\includegraphics[width=6.75cm]{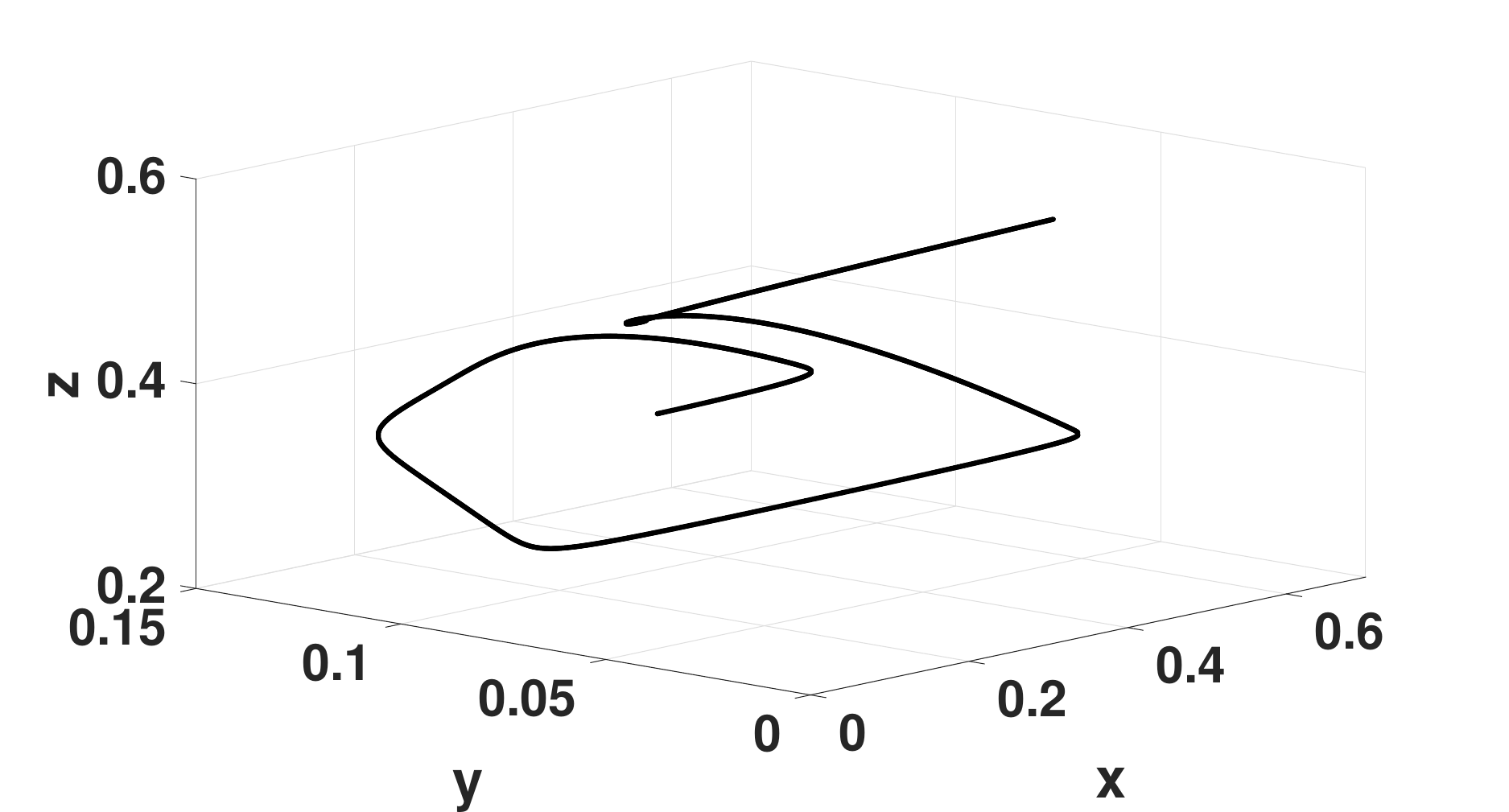}}
  \caption{Bistable dynamics exhibited by system (\ref{nondim3}) at $\beta_1=0.1923$ and $h=0.2648$. The system either approaches a periodic orbit (red) born out of torus  bifurcation  in (A) or to the equilibrium $E_{xz}=(0.357, 0, 0.615)$ in (B). The other parameter values are as in (\ref{parvalues1}).  The initial conditions in (A) are $(0.2995, 0.1175, 0.4154)$ and $(0.34279, 0.104, 0.346)$ in (B).  Note the difference in approach to the periodic orbit  in (A) with figure \ref{orbits_near_FH}(A).}
  \label{bistabilty_h_2648}
\end{figure}
 Below the TR curve and above the yellow curve in figure \ref{two_par_beta2_6}, the equilibrium $E_{xz}$ is locally asymptotically stable and the system exhibits bistable dynamics with periodic orbits born out of torus bifurcation as shown in figure \ref{bistabilty_h_2648}. The coexistence equilibrium $E_2^*$ in this regime is saddle-focus, with one dimensional unstable manifold, $W^u(E_2^*)$, and two-dimensional stable manifold $W^s(E_2^*)$. The approach of a trajectory to the periodic orbit that starts in a neighborhood of $E_2^*$  is shown in  figure  \ref{bistabilty_h_2648}(A). 
Note the difference in  the approaches to the periodic orbits between figures \ref{orbits_near_FH}(A) and \ref{bistabilty_h_2648}(A).

In the region bounded by the PD curve in figure \ref{two_par_beta2_6},  aperiodic and periodic MMO orbits of different signatures  are observed.  Interesting new MMO patterns, not of the standard types  seen in systems with two slow variables and one fast variable, featuring significant variations (by orders of magnitude) in the number and amplitude of  SAOs  between successive spikes (see figure \ref{mmo_h_1825}) are observed. These complicated dynamics demonstrate the uncertainty in the return time between population outbreaks, a realistic phenomena seen in many natural populations. These dynamics will be briefly discussed in the next section.

\section{THE GEOMETRIC SINGULAR PERTURBATION  THEORY APPROACH}

  In this section, we will take a geometrical singular perturbation approach to analyze the complex dynamics exhibited by system $(\ref{nondim3})$ (or equivalently system $(\ref{nondim2})$). The foundation of such geometric approach to analyze systems with a clear separation in time scales was given by Fenichel \cite{F}. We will also review the appropriate canard theory to explain the mechanism behind the MMO patterns.

As $\zeta \to 0$ the trajectories of (\ref{nondim2}) during fast epochs approach to the solutions of the ``layer equations" given by
\begin{equation}\label{layer}
     \left\{
        \begin{array}{r c l}
            {x'}&= &x\phi(x,y,z)\\
 {y'}&=&0\\
    {z'} &= &0.\\
        \end{array}
        \right.
\end{equation}
On the other hand, during slow epochs  trajectories of (\ref{nondim3}) converge to the solutions of the ``reduced problem" given by
\begin{equation}\label{reduced}
        \left\{
        \begin{array}{r c l}
       0 &= &x \phi(x,y,z)\\
   \dot{y}&= &y\chi(x,z)\\
     \dot{z} &=& z\psi(x,y,z).\\
        \end{array}
        \right.
    \end{equation}
    
The subsystems (\ref{layer}) and (\ref{reduced}) will be used to study the dynamics of the full system (\ref{nondim2}) or (\ref{nondim3}) (see \cite{SCT} for details).  

\subsection{The reduced and layer problems} 
\begin{figure}[h!]     
  \centering 
  \subfloat[]{\includegraphics[width=8.0cm]{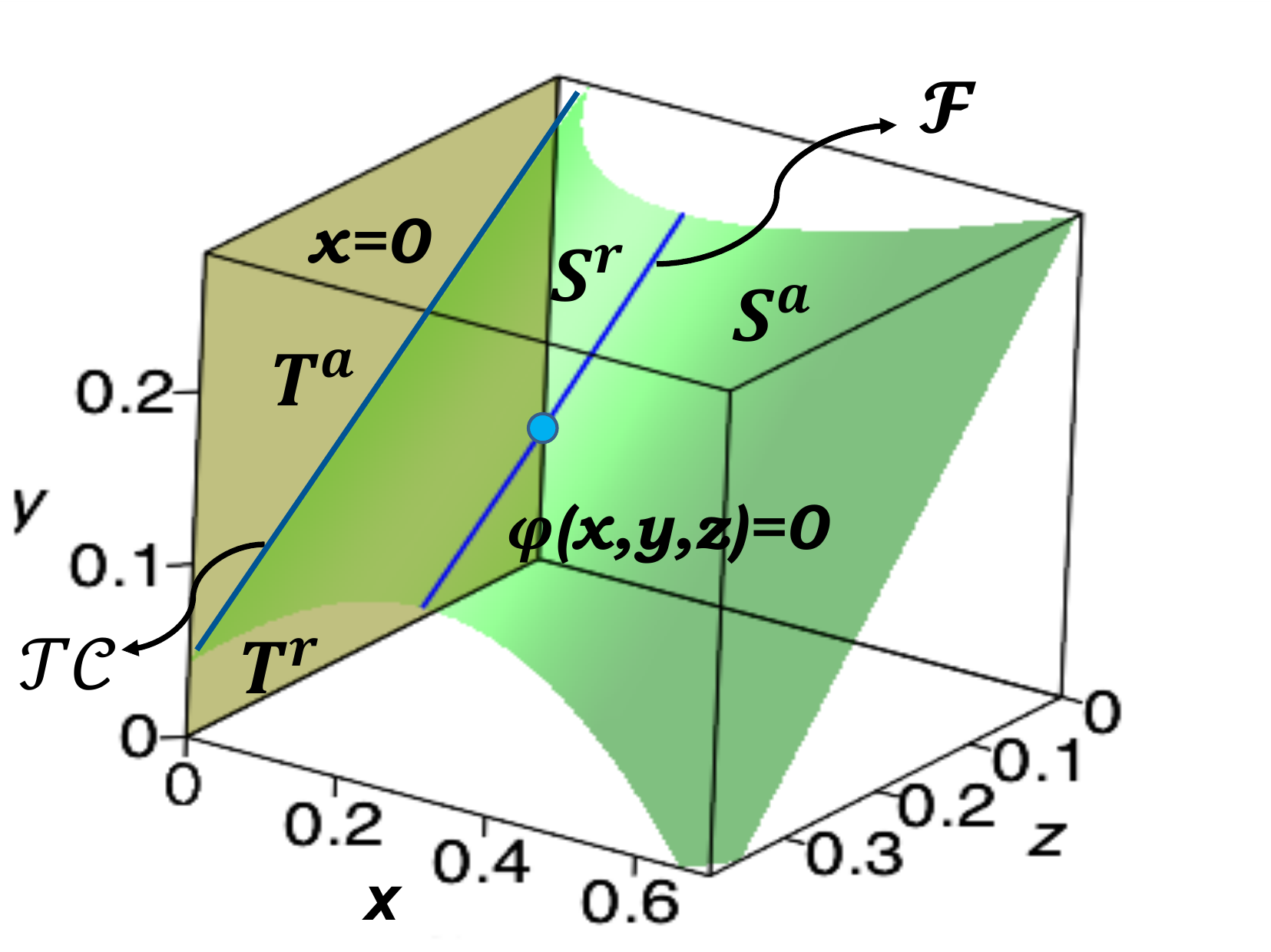}}\quad
  \subfloat[]{\includegraphics[width=7.5cm]{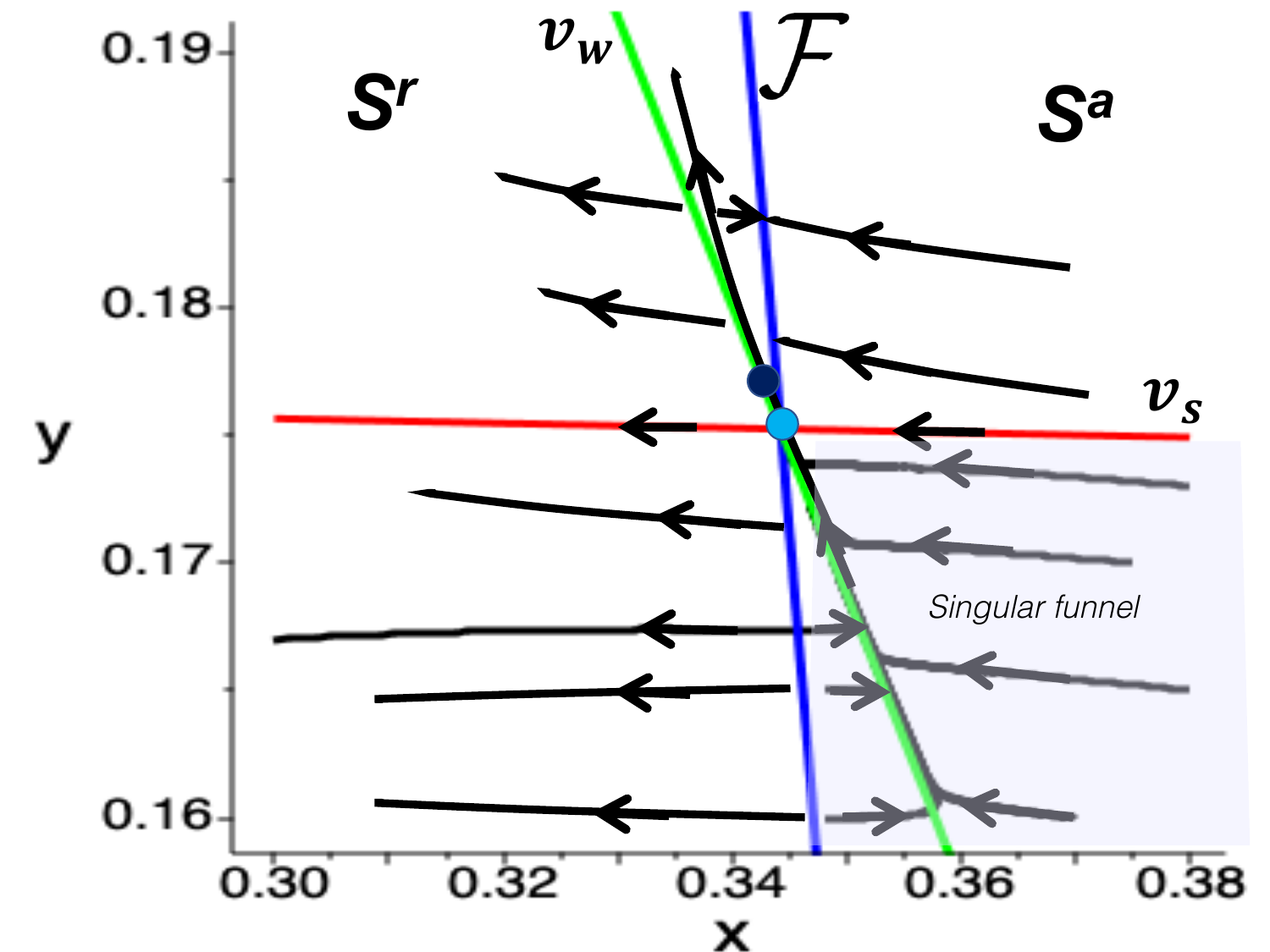}}
 \caption{(A): The critical set $\mathcal{M} =T\cup S$. The folded node is marked by a cyan dot. (B): Reduced dynamics of system (\ref{red2}) zoomed near a folded node singularity (cyan)  projected on the manifold $S$. The black curves are trajectories of  system (\ref{red2}). Also shown are the  ordinary singularity (black), weak eigendirection (green), the strong eigendirection (red) and the singular funnel. The parameter values are $\beta_1=0.25$, $\beta_2=0.35$, $c=0.4$, $d=0.21$, $\alpha_{12}=0.5$, $\alpha_{21}=0.1$, $h=0.8$.} 
 \label{crit_manifold_singular_funnel}
\end{figure}

The algebraic equation in (\ref{reduced}) defines the {\em{critical set}}
\bess \mathcal{M}=\left\{(x,y,z): x=0 {\mathrm{~or~}} \phi(x,y,z)=0\right\} := T\cup S, 
\eess
where  $T =\{ (0,y,z): y, z \geq 0\}$ and $S=\{(x,y,z) \in {\mathbb{R}^3}^+:\phi(x,y,z)=0\}$. The critical set $\mathcal{M}$, consisting of two pieces $T$ and $S$, is the nullsurface of the layer system (\ref{layer}).  The two pieces meet along the line $\mathcal{TC}=\{(0, y, z): y/\beta_1+z/\beta_2=1\}$ which divides the plane $T$ into two normally hyperbolic sheets $T^a=\{(0, y, z):\phi(0, y, z)<0\}$ and $T^r=\{(0, y, z): \phi(0, y, z)>0\}$. The surface $S$ is also divided into two normally hyperbolic sheets $S^a=S\cap \{ \phi_x(x, y, z)<0\}$ and $S^r=S\cap \{ \phi_x(x, y, z)>0\}$  by the curve $\mathcal{F} = S\cap \{ \phi_x(x, y, z)=0\}$.  A linearization of (\ref{layer}) at $\mathcal{M}\setminus \mathcal{F}$ and  $\mathcal{M}\setminus \mathcal{TC}$ yield that the sheets $T^a$ and $S^a$ are attracting while $T^r$ and $S^r$ are repelling (see figure \ref{crit_manifold_singular_funnel}(A)). Saddle-node and transcritical bifurcations of equilibria of the fast subsystem (\ref{layer}) occur along the fold curve $\mathcal{F}$  and the transcritical curve $\mathcal{TC}$ respectively.  By Fenichel's theory \cite{F}, the normally hyperbolic segments of the critical manifold $\mathcal{M}$ perturb to locally invariant attracting and repelling slow manifolds ${T_{\zeta}^a}\cup{S_{\zeta}^a}$ and $T_{\zeta}^r\cup S_{\zeta}^r$ respectively for $\zeta> 0$, and the slow flow restricted to these manifolds is an $O(\zeta)$ perturbation of the reduced flow on $\mathcal{M}$. However, the theory breaks down in neighborhoods of $\mathcal{F}$ and $\mathcal{TC}$ and interesting dynamics  such as relaxation oscillations and MMOs occur.

%

Note that the reduced flow is restricted to the plane $T$ or to the surface $S$.  On the plane $T$, the reduced dynamics solves the system
\begin{eqnarray}\label{plane}    \left\{
\begin{array}{ll} x &=0\\
    \dot{y}&=y\chi(0,z)\\
    \dot{z}&=  z\psi(0,y,z).
       \end{array} 
\right. 
\end{eqnarray}
Since $\dot{y} <-cy$ and $\dot{z} <-dz$ on $T$ with $(0, 0, 0)$ being the global attractor of (\ref{plane}), the reduced flow descends along this plane and approaches $(0, 0, 0)$. As the reduced flow descends, it crosses $\mathcal{TC}$  from $T^a$ to $T^r$ with finite speed, giving rise to singular canards.  The plane $T$ is invariant for all $\zeta>0$, hence the connection between stable and unstable persists for $\zeta>0$, i.e. canards persist for the full system.

%

Since $\phi_z \neq 0$, by the implicit function theorem, the surface $S$ can be locally written as a graph of $z=\theta(x,y)$,  i.e. $\phi(x,y, \theta(x,y))=0$. Hence we can project the dynamics of (\ref{reduced}) onto the $(x,y)$ coordinate chart. Differentiating $\phi(x,y,z)=0$ implicitly with respect to time gives us the relationship $\phi_x\dot{x}+\phi_y\dot{y}+\phi_z\dot{z}=0$. Thus, the reduced flow (\ref{reduced}) restricted to $S$, where $S$ is considered as the graph of $z=\theta(x,y)$ reads as
\begin{eqnarray} \label{red2}  \begin{pmatrix}
 -\phi_x\dot{x}  \\
 \dot{y} 
 \end{pmatrix}= \begin{pmatrix}
 \phi_yy \chi+\phi_zz \psi  \\
y\chi
 \end{pmatrix}\bigg|_{z=\theta(x,y)}.
\end{eqnarray}
System (\ref{red2}) has singularities when $\phi_x=0$ and its solutions blow-up in finite time at $\mathcal{F}$. Hence standard existence and uniqueness results do not hold. Points on $\mathcal{F}$ for which $\phi_yy \chi+\phi_zz \psi\neq 0$ are called ``jump points" and they satisfy the ``normal switching condition" \cite{BKW}. At these points, a solution of (\ref{nondim3}) exits into relaxation after reaching  $\mathcal{F}$ giving rise to {\em{relaxation dynamics}} or more commonly, referred to as  boom and bust cycles in ecology  \cite{SCT} . 
On the other hand,  points of $\mathcal{F}$ where the normal switching condition is violated can give rise to {\em{canards}} as discussed below.  

To analyze the solutions where the normal switching condition fails, we rescale the time $s$ by a phase-space-dependent time transformation factor $-\phi_x$, i.e. $ds = -\phi_x dt_s$ \cite{DGKKOW}. This removes the finite-time blow up of solutions and system (\ref{red2}) transforms to the desingularized system
\begin{eqnarray} \label{desing}  \begin{pmatrix}
 \dot{x}  \\
 \dot{y} 
 \end{pmatrix}= \begin{pmatrix}
\phi_yy \chi+\phi_zz \psi  \\
-\phi_xy\chi
 \end{pmatrix}\bigg|_{z=\theta(x,y)},
\end{eqnarray}
where the overdot denotes $t_s$ derivatives. System (\ref{desing}) is topologically equivalent to system (\ref{red2}) on $S^a$. However, the phase-space-dependent time transformation reverses the orientation of the orbits on $S^r$, therefore, the flow of (\ref{red2}) on $S^r$ is obtained by reversing the direction of orbits of (\ref{desing}). Hence  the reduced flow is either directed towards the fold or away from it (see figure \ref{crit_manifold_singular_funnel}(B)).  

\begin{figure}[h!]     
  \centering 
  {\includegraphics[width=8.75cm]{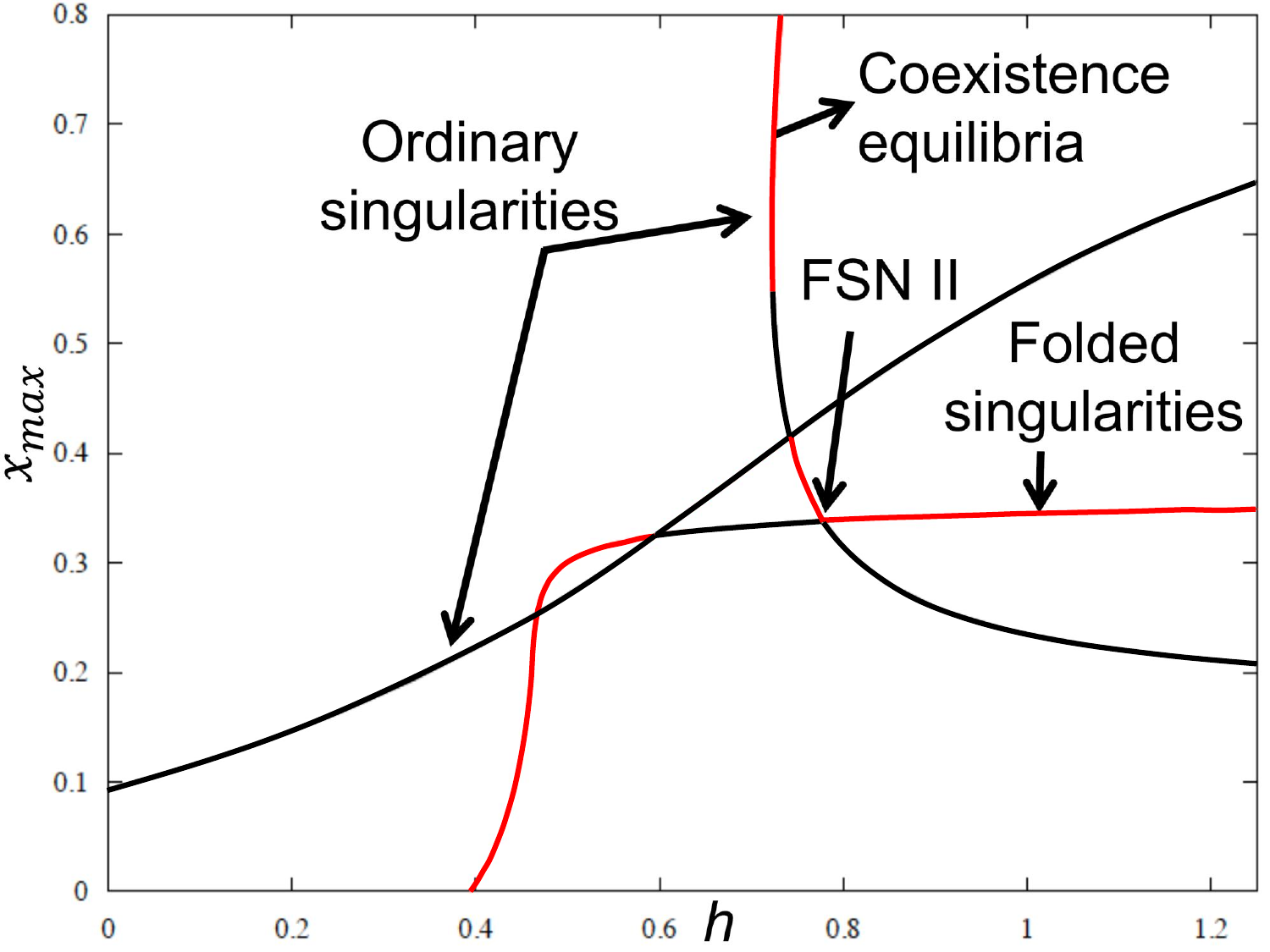}}
 \caption{One-parameter bifurcation  of the desingularized system (\ref{desing}) with varying $h$. The $y$-axis represents the maximum of $x$. The stable parts of a branch are represented in red, while the unstable parts in black. The other parameter values are as in figure \ref{crit_manifold_singular_funnel}. }
 \label{one_par_desing}
\end{figure}

The set of equilibrium points of (\ref{desing}) that do not lie on the fold curve $\mathcal{F}$ (i.e. for which $\phi_x \neq 0$) are {\em{ordinary singularities}}. Such points lie on the curves $\{\chi=0\}\cap \{z=0\}\cap \{\phi_x \neq 0\}$, $\{\psi=0\}\cap \{y=0\} \cap \{\phi_x \neq 0\}$ and $\{\chi=0\}\cap \{\psi=0\}\cap \{\phi_x \neq 0\}$ (see figure \ref{one_par_desing}). On the other hand, equilibrium points of (\ref{desing}) that lie on  $\mathcal{F}$ are  {\em{folded singularities}} or canard points \cite{DGKKOW, SW}. The set of folded singularities form isolated points of $\mathcal{F}$. The trajectories of the reduced system (\ref{reduced}) may pass through the canard points and can thus  cross from $S^a$ to $S^r$ with finite speed, giving rise to singular canards. The classification of a folded singularity as a folded node or a folded saddle or a folded focus or a degenerate folded node is based on the linearization of the folded singularity when considered as an equilibrium of (\ref{desing}) \cite{SW}. More precisely, if $\lambda_{1,2}$ are eigenvalues of the linearization of (\ref{desing}) at a  folded singularity $p_f$, then $p_f$ is a folded focus if $\lambda_{1,2} \in \mathbb{C}$. If $\lambda_{1,2} \in \mathbb{R}$ such that $\lambda_{1}\lambda_2<0$, then $p_f$ is a folded saddle, while if $\lambda_{1}, \lambda_2<0$, then $p_f$ is a folded node. Further degeneracies may occur if one of the eigenvalues pass through $0$ and can give rise to folded saddle node (FSN) bifurcation of types I and II  \cite{DGKKOW}. 

Figure \ref{one_par_desing} shows the curves of ordinary and folded singularities of system (\ref{desing}) as a function of the control parameter $h$, where the other parameter values are held constant. There exists two curves of ordinary singularities, namely the equilibria curves with $y$ absent and the curve of positive  equilibria (coexistence equilibrium state). The latter intersects with the curve of folded singularities at a transcritical bifurcation, also referred to as the FSN II bifurcation. Here the equilibrium $E^*$ of the full system (\ref{nondim2}) and a folded singularity of the desingularized system (\ref{desing}) merge together and then split again, interchanging their type and stability. The equilibrium $E^*$ crosses the fold curve $\mathcal{F}$ at this bifurcation.  The folded singularity switches to a folded node from a folded saddle, and $E^*$ switches from ordinary node to ordinary saddle.

The folded node possesses strong and weak eigenvalues $\lambda_s$ and $\lambda_w$, respectively. The singular strong canard $\gamma_s$ is the unique trajectory to (\ref{desing}) tangent to the strong eigendirection $v_s$, at $p$, while the singular weak canard $\gamma_w$ is tangential to the weak eigendirection, $v_w$. The fold curve $\mathcal{F}$ and the strong singular canard $\gamma_s$ form a trapping region ({\emph{singular funnel}}) on $S^a$ as shown in figure \ref{crit_manifold_singular_funnel}(B), such that all solutions in the funnel converge to the folded node $p$. In fact, the folded node allows a sector of singular canards to flow from $S^a$ to $S^r$ and can give rise to local oscillations for $0<\zeta <<1$  \cite{DGKKOW}. 

 \subsection{Singular orbits, MMOs and relaxation oscillations}The dynamics of the layer problem (\ref{layer}) and the reduced flow (\ref{reduced}) can be combined to construct singular orbits. Such orbits are used to study oscillatory  dynamics such as relaxation oscillations or MMOs of the full system (\ref{nondim3}) \cite{K11}. Typically, a singular periodic orbit for relaxation oscillations in slow-fast systems with ``S-shaped" critical manifolds with two folds is constructed by concatenating the reduced flow, occurring along the attractive branches of the critical manifold, with fast fibers, along which transitions between the branches occur when jump points on the fold curves are reached.  However, the existence of singular canards in $T$ makes the construction of the singular orbit more subtle for system (\ref{nondim3}).  In this case, the reduced flow on $S$ occurs along $S^a$ until it reaches a jump point on $\mathcal{F}$, where it gets connected to $T^a$ by a   fast fiber of (\ref{layer}). On  $T$, as discussed before, the reduced flow  crosses $\mathcal{TC}$ as it descends along $T^a$ and stays on $T^r$ for a while, until it reaches a point $ (0, y_{\tau_0}, z_{\tau_0}) \in T^r$, where a fast orbit concatenates with it. This phenomenon of delay is referred to as the Pontryagin's delay of stability loss \cite{AS, P, MR1, S}. The delay map $\mathcal{P}^0: T^a \to T^r$ is defined by $\mathcal{P}^0(y_0,z_0) = (y(\tau_0(y_0, z_0)), z(\tau_0(y_0, z_0)))$, where the delay $\tau_0$ is expressed by the integral
\bess
\int_0^{\tau_0} \phi(0, y(s), z(s))\ ds =0,
\eess
where $(y(s), z(s))$ solves (\ref{plane}) with initial value $(y_0, z_0) \in T^{a}$. 
 The existence of the Pontryagin's delay of stability loss point on $T^r$ is incorporated in construction of the slow piece of the singular orbit \cite{AS, S}. Singular relaxation cycles are then constructed as continuous concatenation of layer and reduced flow orbit segments and are singular representation of relaxation oscillation cycles. 
   In a similar way, singular orbits which pass through the singular funnel and filtered through folded nodes are singular representation of MMO orbits \cite{BKW}. 

In the parameter regime, where folded node singularities exist for the desingularized system (\ref{desing}), complex oscillatory dynamics such as mixed-mode oscillations (MMOs) are observed. 
The SAOs associated with the MMOs in this model occur near the fold $\mathcal{F}$. System (\ref{nondim3}) exhibits MMOs for different combinations of the parameter values as seen in figures \ref{timeseries_example}, \ref{timeseries_08_819} and \ref{orbits_near_FH_MMO}. 
When $\zeta>0$, a singular Hopf bifurcation occurs at a distance $O(\zeta)$ in parameter space from FSN II bifurcation (see figure \ref{two_par}). Moreover, the desingularized system (\ref{desing}) possesses a folded node singularity (see figure \ref{one_par_desing}). As discussed earlier, according to the Fenichel's theory, the attracting and repelling sheets of the critical set $\mathcal{M}$, perturb smoothly to locally invariant manifolds $S^a_{\zeta}$, $S^r_{\zeta}$, $T^a_{\zeta}$ and $T^r_{\zeta}$ away from the fold curve $\mathcal{F}$ and the transcritical curve $\mathcal{TC}$ for sufficiently small $\zeta>0$. Extending these perturbed slow manifolds by the flow into the vicinity of a folded node singularity  leads to locally twisted intersection of the attracting and the repelling sheets \cite{DGKKOW}, generating the small amplitude oscillations in an MMO orbit. 
A trajectory originating in $S^a_{\zeta}$ to one side of the strong canard is trapped by $S^r_{\zeta}$ and returns towards $S^a_{\zeta}$. This leads to small rotations and the  {\emph{primary weak canard}} serves as the axis of rotation.  A trajectory that originates in  $S^a_{\zeta}$ to the other side of the strong canard is prevented by $S^r_{\zeta}$ from returning  to $S^a_{\zeta}$ and follows a fast direction as it jumps to $T^a_{\zeta}$. This region of rotation is the funnel. Figure \ref{crit_manifold_singular_funnel}(B) gives a view of the  funnel in the singular limit as $\zeta \to 0$. The inset in figure \ref{timeseries_08_819}(A)  gives a zoomed view of the local dynamics near the folded node singularity for $\zeta>0$.  As a trajectory gets trapped in the funnel, it initially experiences rotations imparted by the {\emph{primary weak canard}}. This behavior persists during its passage through the funnel until it exits the vicinity of the folded node. Such small amplitude oscillations are referred to as canard-induced oscillations \cite{BKW, CR}. As the trajectory exits the funnel, the saddle-focus equilibrium $E_1^*$ which lies in $S^r$ plays a role in generating additional oscillations. The local vector field around the equilibrium $E_1^*$ and more precisely, the unstable manifold $W^u(E_1^*)$ now organizes the SAOs. These oscillations are referred to as singular Hopf induced SAOs.  For MMOs to persist,  a trajectory that eventually exits the funnel must return within the region of attraction inside the funnel. In other words, there must exist a global return mechanism. A singular periodic orbit that starts with a fast fiber segment at the folded node and returns within the singular funnel of the folded node formed from the concatenation of orbits of the reduced and layer problems is a candidate for the global return map. For parameter values considered in figure \ref{timeseries_08_819}, such a singular orbit exists.

It is also worth mentioning that  the SAOs associated with the MMO orbit in figure \ref{timeseries_example}  appear quite different from figure \ref{timeseries_08_819}. The parameter regime chosen for both figures is just past the cascade of period doubling bifurcation of the small amplitude limit cycle $\Gamma_h$. However the structure  of the slow manifold $S^{\zeta}$ changes with $\beta_1$, and the amplitude of  SAOs are governed by the extent to which $S_a^{\zeta}$ and $S_r^{\zeta}$ twist while guiding the path of a trajectory through the folded node regime. A tight twisting of $S_a^{\zeta}$ and $S_r^{\zeta}$ forces a decrease in amplitude of the SAOs as observed in figure \ref{timeseries_08_819}.
\begin{figure}[h!]     
  \centering 
  \subfloat[]{\includegraphics[width=7.75cm]{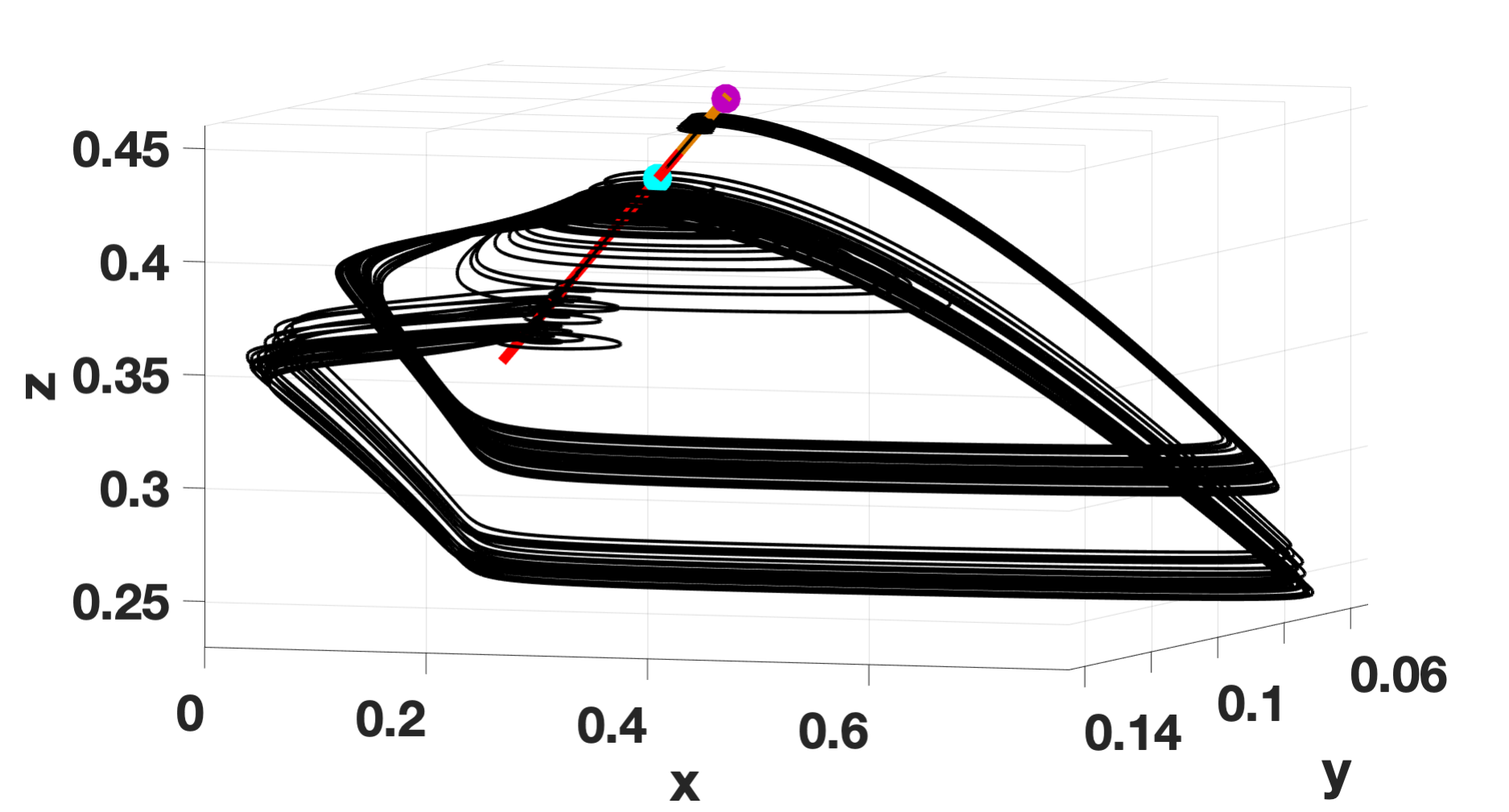}}\qquad
\subfloat[] {\includegraphics[width=6.75cm]{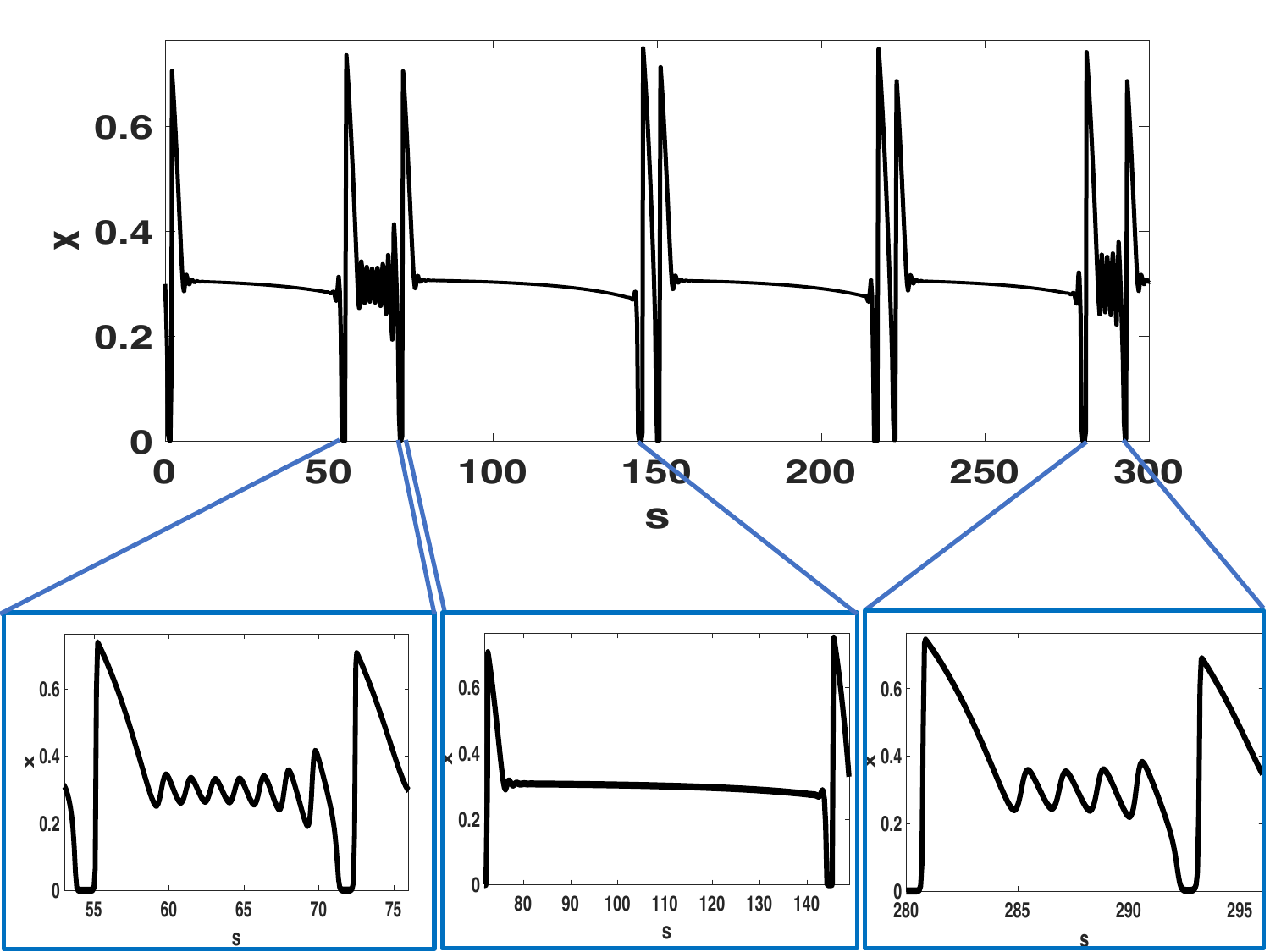}}
  \caption{(A): A chaotic MMO orbit rotating around the weak eigendirection (red) after leaving the equilibrium (magenta) along its unstable eigendirection. Also shown is  the folded node in cyan. (B): Corresponding time series in $x$. Note the amplitudes and epochs of small oscillations  between two large oscillations vary by orders of magnitude. Here $\beta_1=0.1825$, $h =0.2648$  and the other parameter values are as in (\ref{parvalues1}). }
  \label{mmo_h_1825}
\end{figure}

More interestingly, for parameter values within the PD curve  in figure \ref{two_par_beta2_6}, chaotic MMOs with two prominent epochs of SAOs are observed as shown in figure \ref{mmo_h_1825}. The number and amplitude of the small oscillations of the MMOs could vary significantly from one LAO cycle to the next and have a combined  characteristics of the SAOs seen in figure \ref{timeseries_example}, figure \ref{timeseries_08_819} and figure \ref{bistable_attract}(C).  The exact mechanism behind these MMO dynamics is left for future study.

\subsection{Transient MMO dynamics revisited} In Section 3, we showed that prolonged transients in form of MMOs are observed in an immediate neighborhood of FSN II bifurcation. These transients lasted for thousands of generations before the system reached its asymptotic state, namely a periodic orbit born out of a supercritical Hopf bifurcation. We noted that  the SAOs in the transient MMO dynamics have long epochs, where some of the oscillations are too small to be detectable (see figure \ref{bistable_attract}(C)). These MMOs are aperiodic and their chaotic nature can be illustrated by considering a suitable Poincar\'e map near the folded node singularity. During its passage through the funnel, these orbits visit different rotational sectors defined by the secondary canards \cite{DGKKOW}, generating MMOs of mixed signatures $1^{i_1} 1^{i_2}\ldots 1^{i_k}$, where $i_i \in \mathbb{N}$ is very large. To this end, we consider  the Poincar\'e section $\Sigma_h$, defined  by a plane transverse to the critical manifold $S$  through the folded node $p_n$ containing the weak eigenvector.
 \begin{figure}[h!]     
  \centering 
{\includegraphics[width=10.95cm]{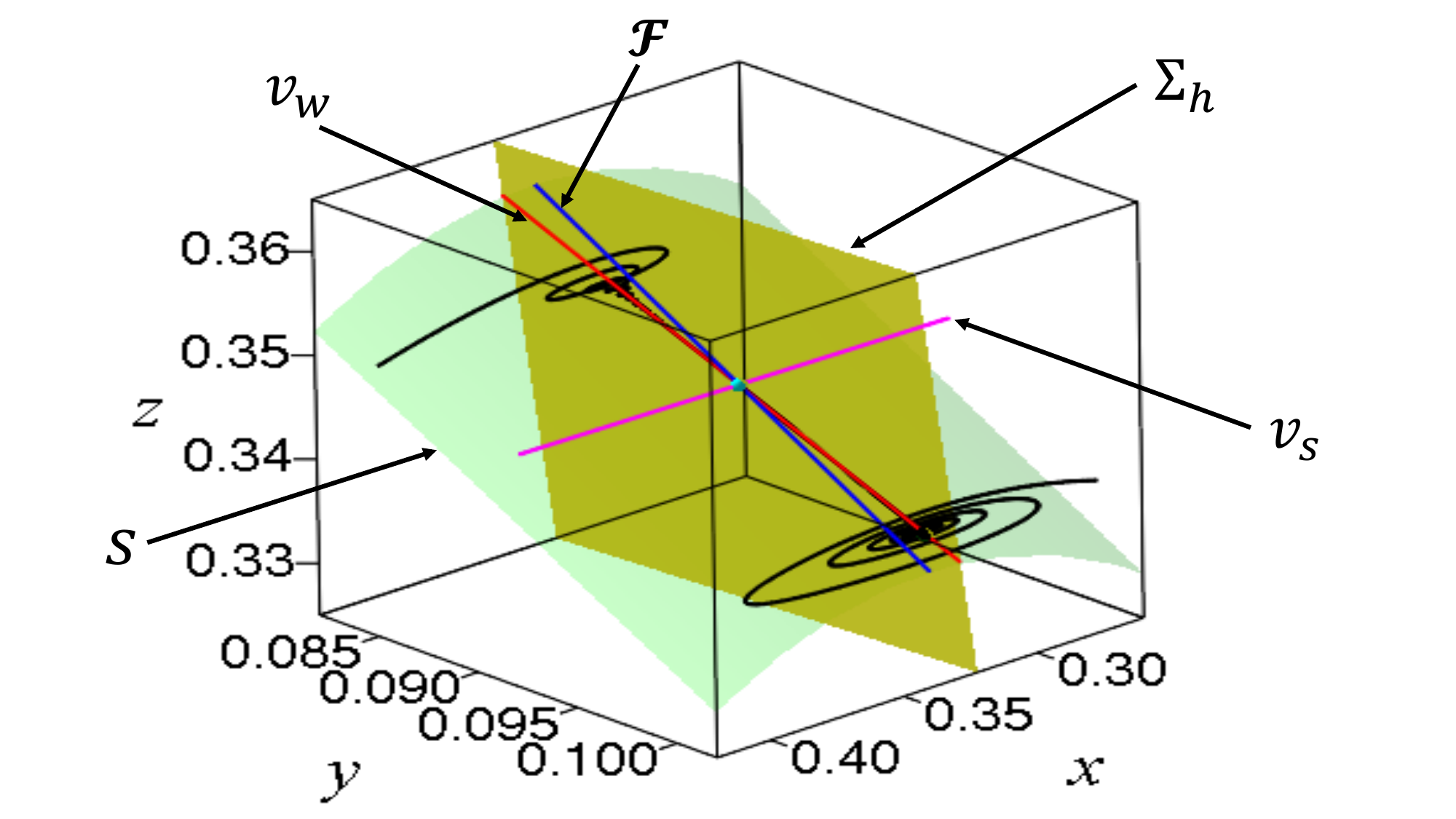}}
\caption{The Poincar\'e section $\Sigma_h$ drawn transverse to the manifold $S$. Also, shown are the strong and weak eigendirections $v_s$ and $v_w$ respectively, along with folded node singularity in cyan. Here $\beta_1=0.25$, $h=0.785$  and other parameter values as in (\ref{parvalues}).}
 \label{poincare_plane}
\end{figure}

 \begin{figure}[h!]     
  \centering 
      \subfloat[Poincar\'e map]{\includegraphics[width=7.05cm]{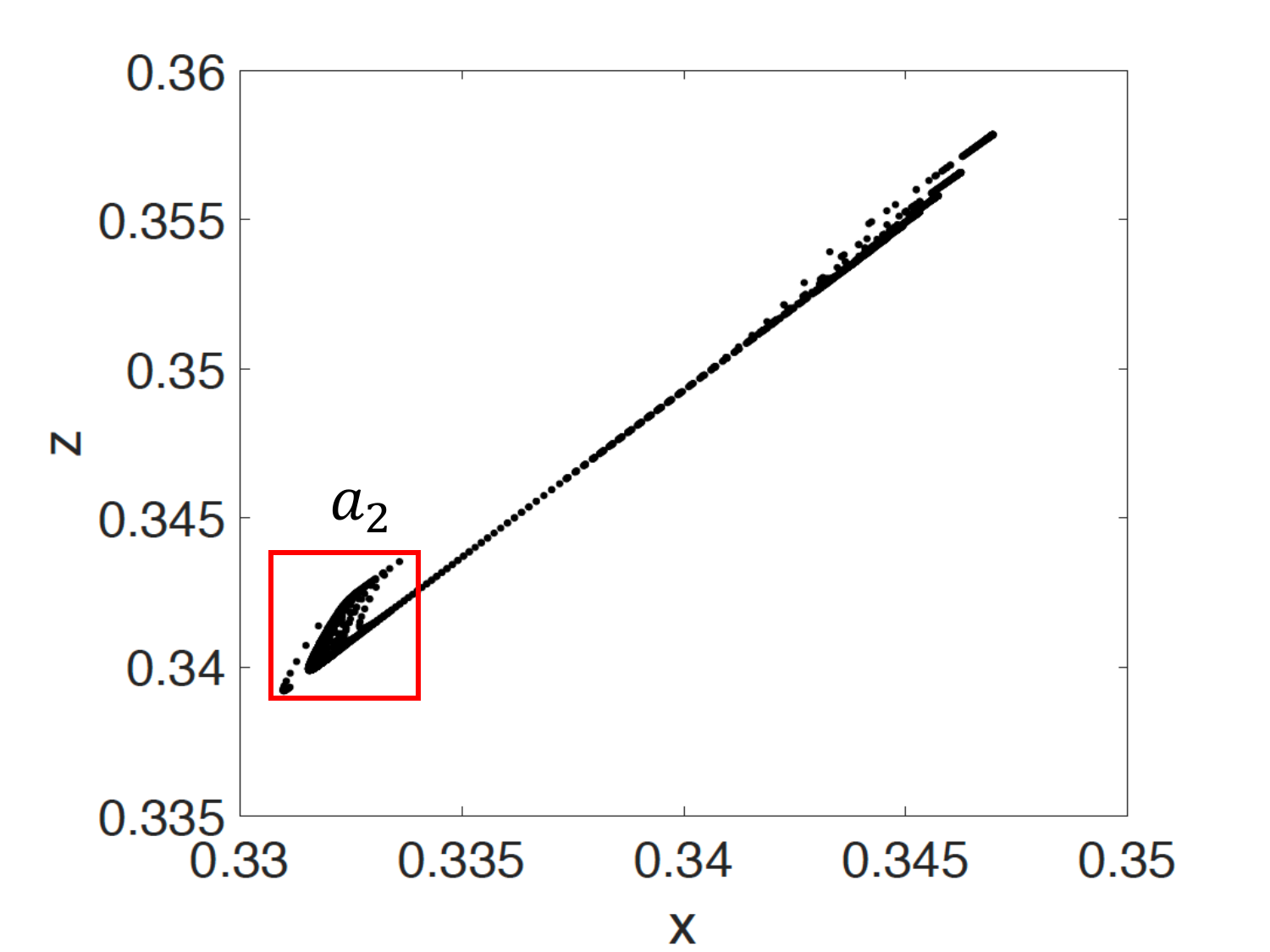}}  \hspace{0.5in}
              \subfloat[Zoomed view of region $a_2$.]{\includegraphics[width=7.05cm]{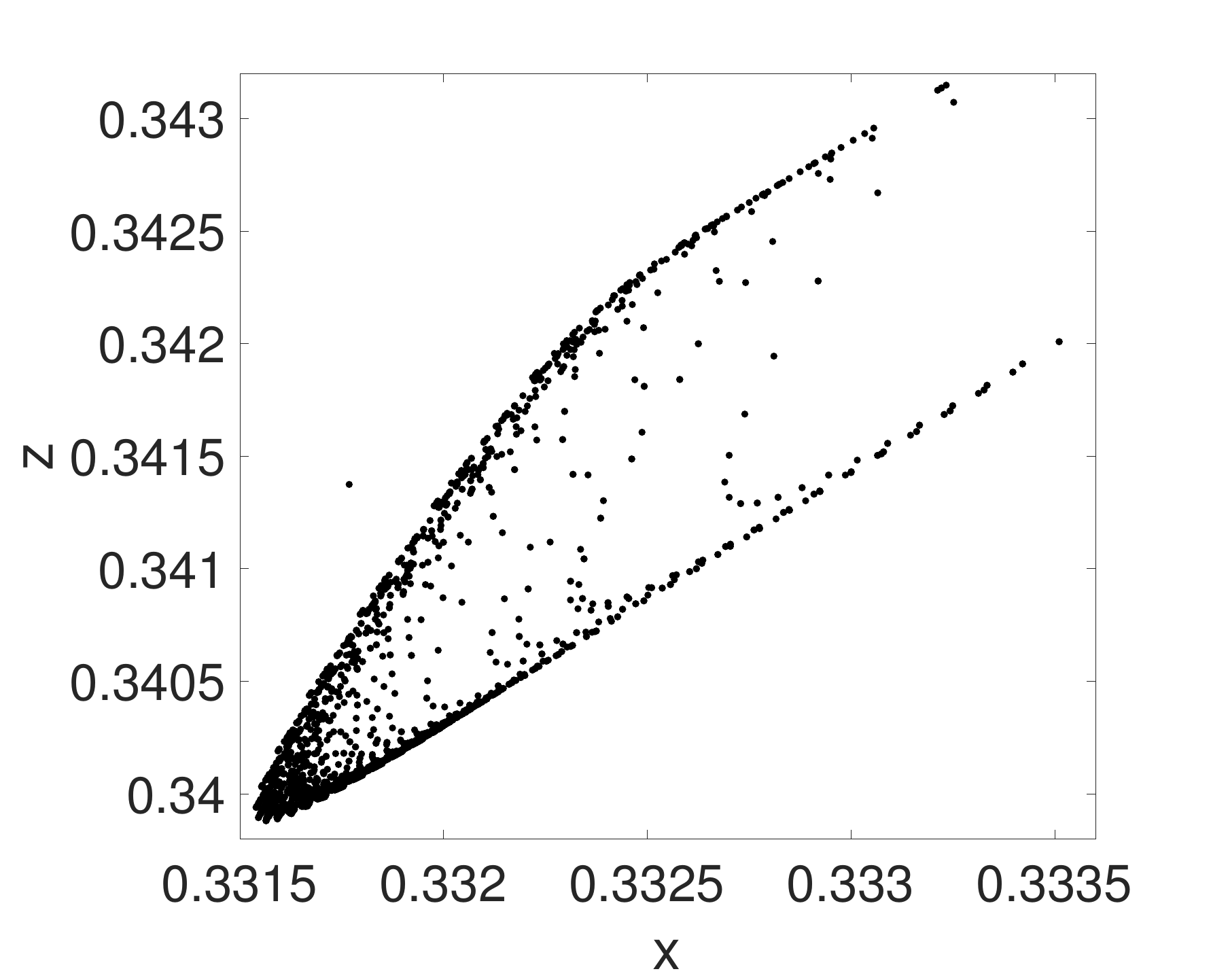}}
  \caption{(A)-(B): $(x, z)$ coordinates of intersection of a trajectory starting at $(0.01, 0.01, 0.12)$  with the cross-section $\Sigma_h:=2.94x+ 1.3y -1.52z = 0.588$ such that $dx/ds<0$.  A zoomed view of the region $a_2$ near the folded node. }
  \label{poincare_sect_sao}
\end{figure}

   A zoomed view of the SAOs of an MMO trajectory along with $\Sigma_h$ is included in figure \ref{poincare_plane}. During the epoch of SAOs, the $y$-coordinate increases while the $z$-coordinate decreases as the trajectory moves toward to $E_1^*$, until the trajectory gets repelled by $W^u(E_1^*)$ and reaches a jump point on  $\mathcal{F}$. Starting with initial condition  $(0.01, 0.01, 0.12)$, the  $(x, z)$ coordinates of the intersections of this trajectory in the direction of decreasing $x$ with the plane $\Sigma_h$ are recorded. The corresponding  Poincar\'e map  is shown in figure \ref{poincare_sect_sao}. The Poincar\'e map indicates that the trajectory repeatedly enters the funnel and stays in a close neighborhood of the folded node singularity $p_n=(0.3383, 0.0923, 0.3474)$, while twisting
 along the primary weak canard and getting pushed towards the equilibrium $E_1^*=(0.3299, 0.1004, 0 .3378)$.  Moreover, a zoomed view of the  map  shows that the trajectory while moving towards $E_1^*$, intersects with the Poincar\'e section at innumerous points forming a dense array of segments as shown in  figure \ref{poincare_sect_sao}(B). A detailed analysis of the underlying mechanism responsible for transition from the transient state to the asymptotic state is done in the companion  paper \cite{Sadhunew}.

\section{Discussion and future outlook}

In this paper, we used singular perturbation theory, theory of canards and numerical simulations to study complex oscillatory dynamics that arise in a three species predator-prey model in a two-trophic ecosystem.  Taking into account that each species operate on different timescales, separation of timescales is introduced in the model, which gives rise to a slow-fast system.  The presence of intraspecific and direct interference competition terms  and the nonlinear functional responses of the predators contribute to the complexities in the model and make the dynamics more realistic. 
 Two-parameter bifurcation diagrams of the system reveal several local and global bifurcations such as period-doubling, torus, saddle-node bifurcation for equilibria, fold-Hopf and generalized Hopf bifurcations, reflecting the rich structure possessed by the system.

The model studied in this paper is an extension of the model considered in \cite{SCT}, where interesting transient dynamics were observed. One such example was a sudden transition from small amplitude oscillations to large amplitude relaxation oscillations.  Such transients can offer an alternative explanation to regime shifts \cite{Hastings} and appear in many ecosystems \cite{morozov}. In this paper, long lasting transients in form of chaotic MMOs are observed near FSN II bifurcation, which eventually asymptote to a stable limit cycle born out of a supercritical Hopf bifurcation. Though a supercritical  Hopf bifurcation is considered as a``safe bifurcation'' in  the asymptotic  limit, the dynamics exhibited by this model shows that this may not be the case on an ecological timescale.   
The presence of a folded node singularity in vicinity of an unstable equilibrium of saddle-focus type gives rise to complexities. The folded node allows certain trajectories on the attracting slow manifold to cross into the repelling slow manifold by creating a funnel. The passage through the funnel induces SAOs in MMOs. Moreover, the proximity of the stable manifold of the saddle-focus equilibrium to the primary weak canard brings a trajectory close to a neighborhood of the equilibrium. As a result, trajectories follow the unstable manifold of the equilibrium leading to additional small rotations before jumping to the other attracting branch of the slow manifold leading to a large excursion in phase space. Thus, finding early warning signs that could help predict when a population will drastically increase is very important.  This requires a detailed understanding of the local mechanism that generates the MMOs and finding the precise locations of the ``jump points" on the fold curve that allow switching from one attracting branch of the slow manifold to the other. This has been considered in a companion paper \cite{Sadhunew}.

The predation efficiencies, $\beta_1, \beta_2$, defined as the ratios of the semi-saturation constants to the carrying capacity of the prey, also played important roles in governing the dynamics. Additional co-dimension 2 bifurcations such as the generalized Hopf bifurcation occurs when one of the predators is considered to be more efficient than the other.  Furthermore,  complicated MMO dynamics with different characteristics timescales and epochs of SAOs between two LAOs are observed as the intraspecific competition strength is varied when  $\beta_2$ is assumed to be significantly higher than  $\beta_1$. These types of oscillations are unique in the framework of MMOs and represent the inherent uncertainties associated with population dynamics and warrants further analysis. The structure of the slow manifold varies with the predation efficiencies and thereby influence the resulting dynamics. In fact, the occurrence (or absence) of MMOs, the epochs of the SAOs in an MMO orbit, and hence formation of canards depend on the twisting properties of the attracting and repelling branches of the slow manifold. Holding the intraspecific competition rate and $\beta_2$ constant, we also noted that as the  predation efficiency of $y$ increases, the competition between the predators go up and eventually leads to competitive exclusion of $y$. The emergence of this boundary equilibrium through torus-like oscillations is remarkable. Such form of oscillations can be related to permanence or persistence in ecological communities and perhaps bear some connection with bursting phenomena, most commonly seen in neuroscience models \cite{hilker}.

 The transient basins, defined as the set of initial conditions which produce chaotic MMO transients whose duration is within a given interval, form complicated structures in phase space, and it is plausible that a small perturbation (such as an environmental fluctuation) can very easily lead to a transition to another basin.  Ecologically, this may be helpful to control population outbreaks, as an external intervention can shift the dynamics from an excited state to a more stable state.

The Koper model \cite{K}, a reduced form of the Hodgkin-Huxley model \cite{DGKKOW}  and  a reduced neuronal competition model \cite{CR} are amongst the few examples of slow-fast systems where a detailed analysis has been carried out near the FSN II points. In \cite{CR}, the authors studied the interaction of canard and singular Hopf mechanism for generation of SAOs near the FSN II bifurcation.  In their work, they classified  the SAOs in an MMO orbit as ``canard-induced" or induced by the local vector field around a saddle-focus equilibrium by numerically computing  the ``way-in/way-out'' function \cite{KW} which describes the maximal delay expected for generic solutions passing through FSN singularity. Similar analysis can be possibly  carried out in this model in the regime far from the fold-Hopf bifurcation, where  the SAOs in an MMO trajectory are distinctly detectable.   Further analysis has been done on the Koper model to study the existence of a homoclinic orbit \cite{GL, MKO}. 
The detection of trajectories with extraordinary large number of small oscillations between large amplitude oscillations and a chaotic return map suggests that the trajectories possibly lie in a small neighborhood of a Shil'nikov homoclinic orbit in this model.  
Explicitly detecting a homoclinic orbit to the saddle-focus equilibrium remains for future study.  The unfolding of the fold-Hopf bifurcation and numerical investigation of global bifurcations responsible for such chaotic dynamics are also subjects of future study. 




\section{Acknowledgement}   I am grateful to Andrew Morozov  for providing constructive feedback on the manuscript. I also thank Chris Cosner and  Sebastian Wieczorek for helpful conversations  at the workshop in ``New Mathematical Methods for Complex Systems in Ecology" at Banff International Research Station in July 2019.  I am thankful to Ryan Lawson, a former student of Georgia College for helping with a numerical code, Saikat Chakraborty Thakur for several helpful discussions and Georgia College Faculty Scholarship Support Program for supporting this research.

\begin{remark}\label{rmk2} Most of the numerical simulations in this paper were done in MATLAB. We used the predefined routine ODE$45$  with relative and absolute error tolerances $10^{-11}$ and $10^{-12}$ respectively.

\end{remark}




\end{document}